\documentclass[11pt]{amsart}

\usepackage{amsmath}
\usepackage{amscd}
\usepackage{amssymb}
\usepackage{latexsym}
\usepackage{stmaryrd} 
\usepackage{mathbbol} 

\usepackage{color} 
\definecolor{darkgreen}{cmyk}{1,0,1,.2}
\definecolor{m}{rgb}{0.5,0.1,1}

\usepackage[arrow,curve,matrix,tips,2cell]{xy}
  \SelectTips{eu}{10} \UseTips
  \UseAllTwocells
\expandafter\def\csname r@tocindent4\endcsname{0pt} 

\newtheorem{theorem}{Theorem}[section]
\newtheorem*{theorem*}{Theorem}
\newtheorem{lemma}[theorem]{Lemma}
\newtheorem*{lemma*}{Lemma}
\newtheorem{corollary}[theorem]{Corollary}
\newtheorem{proposition}[theorem]{Proposition}

\newtheorem{remark}[theorem]{Remark}
\newtheorem{definition}[theorem]{Definition}

\newtheorem{subtheorem}{Theorem}[subsection]
\newtheorem{subdefinition}[subtheorem]{Definition}
\newtheorem{subproposition}[subtheorem]{Proposition}
\newtheorem{sublemma}[subtheorem]{Lemma}
\newtheorem{subcorollary}[subtheorem]{Corollary}
\newtheorem{subremark}[subtheorem]{Remark}
\newcommand{\bgl}{\begin{equation}} 
\newcommand{\egl}{\end{equation}}
\newcommand{\bgloz}{\begin{equation*}} 
\newcommand{\egloz}{\end{equation*}}
\newcommand{\bgln}{\begin{eqnarray}} 
\newcommand{\egln}{\end{eqnarray}}
\newcommand{\bglnoz}{\begin{eqnarray*}} 
\newcommand{\eglnoz}{\end{eqnarray*}}
\newcommand{\btheo}{\begin{theorem}}
\newcommand{\etheo}{\end{theorem}}
\newcommand{\bsubtheo}{\begin{subtheorem}}
\newcommand{\esubtheo}{\end{subtheorem}}
\newcommand{\btheooz}{\begin{theorem*}}
\newcommand{\etheooz}{\end{theorem*}}
\newcommand{\blemma}{\begin{lemma}}
\newcommand{\elemma}{\end{lemma}}
\newcommand{\bsublemma}{\begin{sublemma}}
\newcommand{\esublemma}{\end{sublemma}}
\newcommand{\blemmaoz}{\begin{lemma*}}
\newcommand{\elemmaoz}{\end{lemma*}}
\newcommand{\bproof}{\begin{proof}}
\newcommand{\eproof}{\end{proof}}
\newcommand{\bbew}{\begin{beweis}}
\newcommand{\ebew}{\end{beweis}}
\newcommand{\bremark}{\begin{remark}\em}
\newcommand{\eremark}{\end{remark}}
\newcommand{\bsubremark}{\begin{subremark}\em}
\newcommand{\esubremark}{\end{subremark}}
\newcommand{\bdefin}{\begin{definition}}
\newcommand{\edefin}{\end{definition}}
\newcommand{\bsubdefin}{\begin{subdefinition}}
\newcommand{\esubdefin}{\end{subdefinition}}
\newcommand{\bprop}{\begin{proposition}}
\newcommand{\eprop}{\end{proposition}}
\newcommand{\bsubprop}{\begin{subproposition}}
\newcommand{\esubprop}{\end{subproposition}}
\newcommand{\bcor}{\begin{corollary}}
\newcommand{\ecor}{\end{corollary}}
\newcommand{\bsubcor}{\begin{subcorollary}}
\newcommand{\esubcor}{\end{subcorollary}}
\newcommand{\bfa}{\begin{cases}} 
\newcommand{\efa}{\end{cases}}
\newcommand{\cH}{\mathcal H}

\newcommand{\cJ}{\mathcal J}
\newcommand{\cK}{\mathcal K}
\newcommand{\cL}{\mathcal L}

\newcommand{\cO}{\mathcal O}

\newcommand{\cR}{\mathcal R}

\newcommand{\cU}{\mathcal U}

\newcommand{\cX}{\mathcal X}

\def\Cz{\mathbb{C}}

\def\Zz{\mathbb{Z}}

\def\1z{\mathbb{1}}
\newcommand{\fA}{\mathfrak A}

\newcommand{\fI}{\mathfrak I}

\def\xf{{\bf{x}}}
\def\yf{{\bf{y}}}
\def\zf{{\bf{z}}}
\newcommand{\an}[1]{``#1''} 
\newcommand{\ti}{\tilde}

\newcommand{\lori}{\longrightarrow}

\newcommand{\ma}{\mapsto} 
\newcommand\onto{\twoheadrightarrow} 
\newcommand\into{\hookrightarrow} 
\newcommand{\Rarr}{\Rightarrow} 
\newcommand{\Larr}{\Leftarrow} 
\newcommand{\LRarr}{\Leftrightarrow} 

\newcommand{\ve}{\varepsilon}

\def\SEMI{\mbox{$\times\kern-2pt\vrule height5pt width.6pt \kern3pt $}}

\newcommand{\End}{{\rm End}\,}
\newcommand{\Aut}{{\rm Aut}\,}

\newcommand{\img}{{\rm Im\,}}

\newcommand{\id}{{\rm id}}

\newcommand{\stalg}{{}^*\text{-alg}}

\newcommand{\Ind}{\mathrm{ Ind}\,}

\newcommand{\Ad}{{\rm Ad\,}}

\newcommand{\reg}{^\times} 
\newcommand{\lspan}{{\rm span}} 
\newcommand{\clspan}{\overline{\lspan}} 
\newcommand{\abs}[1]{\lvert#1\rvert} 
\newcommand{\norm}[1]{\left\|#1\right\|} 
\newcommand{\defeq}{\mathrel{:=}} 
\newcommand{\eqdef}{\mathrel{=:}} 

\newcommand{\dop}{\text{: }} 
\newcommand{\fa}{\text{ for all }} 
\newcommand{\ilim}{\varinjlim} 
\newcommand{\e}[1]{e_{\left[#1\right]}} 
\newcommand{\rte}{\overset{e}{\rtimes}} 

\newcommand{\Isom}{{\rm Isom}\,}

\newcommand{\lge}{\left\{} 
\newcommand{\rge}{\right\}} 
\newcommand{\lru}{\left(} 
\newcommand{\rru}{\right)} 
\newcommand{\leck}{\left[} 
\newcommand{\reck}{\right]} 
\newcommand{\lsp}{\left\langle} 
\newcommand{\rsp}{\right\rangle} 
\newcommand{\rukl}[1]{\lru #1 \rru} 
\newcommand{\eckl}[1]{\leck #1 \reck} 
\newcommand{\gekl}[1]{\lge #1 \rge} 
\newcommand{\spkl}[1]{\lsp #1 \rsp} 
\newcommand{\menge}[2]{\gekl{ #1 \dop #2 }} 
\DeclareMathOperator{\res}{res}
\newcommand{\tih}{\ti{h}}
\newcommand{\DP}{D_r(P)}
\newcommand{\DiP}{D_r^{(\infty)}(P)}
\newcommand{\taui}{\tau^{(\infty)}}
\newcommand{\pii}{\pi^{(\infty)}}
\newcommand{\E}[2]{E_{{#1}^{-1} \cdot {#2}}}
\newcommand{\iu}{^{(\infty)}}
\begin{document}

\title[K-theory of semigroup C*-algebras]{On the K-theory of the C*-algebra generated by the left regular representation of an Ore semigroup}

\author[J. Cuntz]{Joachim Cuntz}
\address{Joachim Cuntz, Department of Mathematics, Westf{\"a}lische Wilhelms-Universit{\"a}t M{\"u}nster, Einsteinstra{\ss}e 62, 48149 M{\"u}nster, Germany}
\email{cuntz@uni-muenster.de}
\author[S. Echterhoff]{Siegfried Echterhoff}
\address{Siegfried Echterhoff, Department of Mathematics, Westf{\"a}lische Wilhelms-Universi\-t{\"a}t M{\"u}nster, Einsteinstra{\ss}e 62, 48149 M{\"u}nster, Germany}
\email{echters@uni-muenster.de}
\author[X. Li]{Xin Li}
\address{Xin Li, Department of Mathematics, Westf{\"a}lische Wilhelms-Universit{\"a}t M{\"u}nster, Einsteinstra{\ss}e 62, 48149 M{\"u}nster, Germany}
\email{xinli.math@uni-muenster.de}


\thanks{\scriptsize{Research supported by the Deutsche Forschungsgemeinschaft (SFB 878) and by the ERC through AdG 267079.}}

\begin{abstract}
We compute the K-theory of  C*-algebras  generated by the left regular representation of left Ore semigroups satisfying  certain regularity conditions. 
 Our result describes the K-theory of these semigroup C*-algebras in terms of the K-theory for the reduced group C*-algebras of certain groups which are typically easier to handle. Then we apply our result to specific semigroups from algebraic number theory.
\end{abstract}

\maketitle


\setlength{\parindent}{0pt} \setlength{\parskip}{0.5cm}

\section{Introduction}

Let $P$ be a (discrete) semigroup. If $P$ admits left cancellation, then left translation defines an action of $P$ by isometries $V_p$, $p \in P$, on the Hilbert space $\ell^2 P$. When $P$ is a group, the $V_p$ are unitaries and the reduced C*-algebra $C^*_r(P)$ generated by the operators $V_p$ is one of the most classical objects of study in the theory of operator algebras. The analogous C*-algebra for a genuine semigroup has recently, partly triggered by natural examples, found attention and has been studied in various connections. We call them (reduced or regular) semigroup C*-algebras. The interested reader may consult \cite{Li2} for a brief account of the historical background of these C*-algebras attached to semigroups.

The possibility of describing $C^*_r(P)$, for a left cancellative semigroup $P$, as a universal C*-algebra with generators and relations has been analyzed in connection with amenability properties of $P$ in \cite{Li2}. Also, such a description was discussed in detail for the important example of the \an{$ax+b$-semigroup} $R\rtimes R^\times$ for the ring of integers $R$ in a number field in \cite{C-D-L}. In this latter paper also the KMS-structure for a natural one-parameter group on $C^*_r(R\rtimes
R^\times)$ was studied and it was shown that it is partly governed by the ideal class group for $R$.

In the present paper we set out to determine the K-theoretic invariants of $C^*_r(P)$ for a class of semigroups containing the semigroups arising from number theory that we are interested in. Here is our main result:
\btheooz
Let $P$ be a countable left Ore semigroup. Assume that the family of constructible right ideals $\cJ$ of $P$ is independent (\S~\ref{ideal-structure}), and that the enveloping group $G$ of $P$ satisfies the Baum-Connes conjecture with coefficients. 
Let $\mathcal I$ denote the $G$-saturation of $\mathcal J\setminus\{\emptyset\}$ in the power set $\mathcal P(G)$
of $G$. Then the K-theory of the semigroup C*-algebra $C^*_r(P)$ can be described as follows:
$$
K_*(C^*_r(P)) \cong \bigoplus_{[X] \in G\backslash \mathcal I} K_*(C^*_r(G_X)),
$$
where $G_X= \menge{g \in G}{g \cdot X = X}$ denotes the stabilizer of $X\in \mathcal I$ under the action of $G$ on $\mathcal I$.
\etheooz
In fact, we only need the Baum-Connes conjecture with coefficients in two specific $G$-C*-algebras. Moreover, in good situations, it turns out that $C^*_r(P)$ and $\bigoplus_{[X] \in G\backslash \mathcal I} C^*_r(G_X)$ are actually KK-equivalent. We refer the reader to \S~\ref{abstract-concrete} for more explanations and more precise formulations of our result. Let us now explain the basic ideas behind the proof:

As a first step, we need an embedding of $C^*_r(P)$ as  a full corner of a (reduced) crossed product $D\rtimes_rG$ of a commutative C*-algebra $D$ by an enveloping group $G$ for $P$ (see Section~\ref{dilations}). The existence of such a crossed product follows from the left Ore condition on $P$ (see \cite{La}). As a consequence, the K-theory of $C_r^*(P)$ is isomorphic to 
the K-theory of $D\rtimes_rG$.

We then prove a rather general K-theoretic theorem which, in many situations, allows to reduce the computation of $ K_*(D\rtimes_rG)$ to the, often much simpler, computation of the K-theory of C*-algebras associated with certain subgroups of $G$. Our key technical result concerns the following situation. Assume
that $D$ is a commutative C*-algebra generated by a multiplicative family $\menge{e_i}{i \in I}$ of projections, satisfying a certain independence condition, and that a group $G$ acts on $D$ leaving the generating family invariant. We then show under the assumption that $G$ satisfies the Baum-Connes conjecture for the coefficient algebras $D$ and $c_0(I)$ that the computation of the K-theory of the crossed product C*-algebra $D\rtimes_r G$ is equivalent to the computation of the K-theory of the much simpler crossed product $c_0(I)\rtimes_r G$ (see Section~\ref{generalK}). The proof uses techniques that have been developed in connection with the Baum-Connes conjecture in \cite{C-E-O}, \cite{E-L-P-W} and \cite{Mey-Ne}. A combined statement of the relevant results is given in \cite{E-N-O}.  Note that  by \cite{HK} all amenable groups (among many others)  satisfy the Baum-Connes conjecture, so the results apply in particular to our 
 motivating examples $R\rtimes R^\times$.

Now, on the other hand, if a group $G$ acts on $c_0(I)$ where $I$ is a discrete set, then simple imprimitivity considerations show that the crossed product $c_0(I) \rtimes_r G$ is Morita equivalent to a direct sum of the (reduced) group C*-algebras of the stabilizer groups.

In the case of the crossed product $D\rtimes_rG$ connected to the left Ore semigroup $P$, the algebra 
$D$ is generated by the set of projections $\{E_X: X\in \mathcal I\}$ with $E_X$ the orthogonal projection from 
$\ell^2G$ to $\ell^2(X)\subseteq \ell^2(G)$. The  independence condition 
for this set of projections follows from a similar 
independence condition on the set of constructible right ideals $\mathcal J$ in $P$. This gives the 
result of our theorem. Moreover, if $G$ satisfies a certain strong version of the Baum-Connes conjecture
(which again holds, among others, for all amenable groups) we can conclude the stronger result  that $C_r^*(P)$ is KK-equivalent 
to the direct sum $\bigoplus_{[X]\in G\backslash \mathcal I} C_r^*(G_X)$. Note that the $G$-orbits in  $\mathcal I$ and the stabilizers 
$G_X$ are easily determined in specific examples.

Under the same assumptions on our semigroup $P$ as above, there exists a natural diagonal map $C^*_r(P) \to C^*_r(P)\otimes_{\min} C^*_r(P)$. This means that, just as for a group C*-algebra, the K-homology of $C^*_r(P)$ becomes a ring via this diagonal map. The KK-equivalence between $C^*_r(P)$ and the direct sum of the C*-algebras of the stabilizer groups, that we construct in presence  of the strong Baum-Connes conjecture, induces in fact an isomorphism of K-homology rings.

As mentioned above, our motivating examples are the semigroups attached to a Dedekind domain $R$, such as the ring of integers in an algebraic number field, or function field, $K$. For such a ring we consider the multiplicative semigroup $R^\times$, the multiplicative semigroup of principal ideals and the $ax+b$-semigroup $R\rtimes R^\times$ (see Section~\ref{Dede}). These semigroups have obvious enveloping groups $K^\times$, the group of principal fractional ideals and $K\rtimes K^\times$. The set $\mathcal I$ which appears when we apply our theorem can be identified with the set of fractional ideals (for both $R\reg$ and the semigroup of principal ideals), or with the translates of fractional ideals, in $K$, respectively. The stabilizer groups are essentially the group of invertible elements in $R^\times$, trivial or the group of invertible elements in $R\rtimes R^\times$. The orbits in $\mathcal I$ for the action of the enveloping group are labeled by the ideal class group $Cl_K$ in each case. We note that in the case of the multiplicative semigroups, there are natural actions of the class group on the K-theory of the corresponding semigroup C*-algebras.

Finally we turn to a study of specific structural properties of the C*-algebra $C^*_r(R\rtimes R^\times)$ for the ring of integers $R$ in a number field (see \S~\ref{ax+b}). This algebra is of special interest for many reasons. As mentioned above, it has an intriguing KMS-structure, but it also has a unique maximal ideal and the quotient by this ideal gives the ring C*-algebra $\fA[R]$ studied in \cite{Cu-Li}. This ring C*-algebra is purely infinite and simple and can be represented as a crossed product by actions on adele spaces in different ways. In \cite{Cu-Li} we had also determined its K-theory for a first class of number fields. The complete K-theoretic computation is obtained in \cite{Li-Lu}.

Using a criterion from \cite{Pas-Ror} we can now show that $C^*_r(R\rtimes R^\times)$ is purely infinite (though of course not simple) and has the ideal property. These properties are of structural interest for a C*-algebra. Using our K-theory computation and another criterion from \cite{Pas-Ror}, we can show that $C^*_r(R\rtimes R^\times)$ on the other hand does not have real rank zero. The first named author is indebted to C. Pasnicu and G. Gong for drawing his attention to these properties.

\section{Preliminaries}
\label{pre}

\subsection{Semigroups}

A semigroup is a set $P$ together with an associative binary operation (or multiplication) $P \times P \to P$; $(p,q) \ma p q$. We will not consider (non-trivial) topologies on our semigroups, which means that topologically, all our semigroups will be viewed as discrete sets. A unit element in a semigroup $P$ is an element $e$ in $P$ with the property that $ep = pe = p$ for all $p$ in $P$. All the semigroups in this paper are assumed to have unit elements. In addition, since we would like to use KK-theory in \S~\ref{generalK}, all our semigroups in \S~\ref{generalK} are supposed to be countable so that the semigroup C*-algebras will be separable.

Moreover, a semigroup $P$ is called left cancellative if for all $p$, $x$ and $y$ in $P$, $px=py$ implies $x=y$. Similarly, a semigroup $P$ is called right cancellative if for all $p$, $y$ and $y$ in $P$, $yp=yp$ implies $x=y$. A semigroup is called cancellative if it is both left and right cancellative.

\subsection{Ideal structure}
\label{ideal-structure}

A left ideal of a semigroup $P$ is a subset $X$ of $P$ which is invariant under left multiplications, i.e. for every $x$ in $X$ and $p$ in $P$, $px$ lies in $X$ again. Similarly, a right ideal of a semigroup $P$ is a subset $X$ of $P$ which is invariant under right multiplications, i.e. for every $x$ in $X$ and $p$ in $P$, $xp$ lies in $X$ again.

In the analysis of semigroup C*-algebras, a certain family of right ideals plays an important role. It is defined as follows:
\bsubdefin
\label{cJ}
For a semigroup $P$, let $\cJ$ be the smallest family of right ideals of $P$ satisfying
\begin{itemize}
\item $\emptyset, P \in \cJ$ 
\item $\cJ$ is closed under left multiplication and taking pre-images under left multiplication ($X \in \cJ, p \in P \Rarr pX, p^{-1}X \in \cJ$)
\item $\cJ$ is closed under finite intersections ($X,Y \in \cJ \Rarr X \cap Y \in \cJ$).
\end{itemize}
Here we define for every subset $X$ of $P$ and for all $p \in P$:
\bgloz
  pX \defeq \menge{px}{x \in X} \text{ and } p^{-1}X \defeq \menge{q \in P}{pq \in X}.
\egloz
\esubdefin
It follows directly from this definition that $\cJ$ consists of $\emptyset$ and arbitrary finite intersections of right ideals of the form $q_1^{-1} p_1 \dotsm q_n^{-1} p_n P$ for $q_1, \dotsc, q_n, p_1, \dotsc, p_n \in P$. Elements in $\cJ$ are called constructible right ideals of $P$.

We need the following
\bsubdefin
The family $\cJ$ is said to be independent (we also say that the constructible right ideals of $P$ are independent) if for all right ideals $X, X_1, \dotsc, X_n$ in $\cJ$ with $X = \bigcup_{j=1}^n X_j$, we must have $X=X_j$ for some $1 \leq j \leq n$.
\esubdefin
In other words, $\cJ$ is independent if for every right ideal $X$ in $\cJ$, the following holds: Given $X_1, \dotsc, X_n$ in $\cJ$ which are proper subsets of $X$ ($X_j \subsetneq X$ for all $1 \leq j \leq n$), then the union $\bigcup_{j=1}^n X_j$ is again a proper subset of $X$ ($\bigcup_{j=1}^n X_j \subsetneq X$).

This independence condition plays an important role when one tries to describe amenability of semigroups in terms of semigroup C*-algebras (see \cite{Li2}). But as we will see, it will also play a crucial role in our K-theoretic computations.

\subsection{Ore semigroups}
\label{Ore}

Our K-theoretic computations only work for so-called left Ore semigroups.

\bsubdefin
A semigroup is called right reversible if every pair of non-empty left ideals has a non-empty intersection.
\esubdefin

\bsubdefin
A semigroup is said to satisfy the left Ore condition if it is cancellative and right reversible. A semigroup with such properties is called a left Ore semigroup.
\esubdefin

The following result is the reason why the left Ore condition is so useful:
\bsubtheo[Ore, Dubreil]
\label{OD}
A semigroup $P$ can be embedded into a group $G$ such that $G = P^{-1} P = \menge{q^{-1} p}{p,q \in P}$ if and only if $P$ satisfies the left Ore condition. In this case, the group $G$ is determined up to canonical isomorphism by the universal property that every semigroup homomorphism $P \to G'$ from $P$ to a group $G'$ extends uniquely to a group homomorphism $G \to G'$.
\esubtheo
When we write $G = P^{-1} P$ in this theorem, then we are identifying $P$ with its image in $G$ under the embedding of $P$ into $G$.

The reader may consult \cite{Cl-Pr}, Theorem~1.23 or \cite{La}, \S~1.1 for more explanations about this theorem. For a left Ore semigoup $P$, let us call the (up to canonical isomorphism unique) group $G$ which appears in the theorem the enveloping group of $P$. It is also called the group of left quotients (which explains the terminology \an{left Ore semigroup}).

Instead of giving a full proof of this theorem, we now describe an explicit model for the enveloping group in order to illustrate an important idea. Let $P$ be a semigroup. We define a partial order on $P$ by setting $p \leq q :\LRarr q \in Pp$. Here $Pp$ is the left principal left ideal of $P$ generated by $p$, i.e. $Pp = \menge{xp}{x \in P}$. It is straightforward to see that $P$ is right reversible if and only if $P$ is upwards directed with respect to this partial order, which means that for all $p_1,p_2 \in P$, there exists $q \in P$ such that $p_1 \leq q$ and $p_2 \leq q$. If we further assume that $P$ is right cancellative, then $p \leq q$ implies that there exists a unique element $r \in P$ with $q = rp$. We denote this element $r$ by $q p^{-1}$. The observations made so far tell us that given a right reversible, right cancellative semigroup $P$, we can form an inductive system of sets indexed by the elements in $P$ ordered by \an{$\leq$} in the following way:
\begin{itemize}
\item for every $p \in P$, the $p$-th set is given by $P$ itself
\item for every $p, q \in P$ with $p \leq q$, the structure map from the $p$-th set to the $q$-th set is given by left multiplication with $q p^{-1}$: $P \to P; x \ma (qp^{-1})x$.
\end{itemize}
We can then form the set-theoretical inductive limit of this system and endow it with a binary operation so that we again obtain a semigroup. Here are the details: As a first step, we take the (set-theoretical) disjoint union $\bigsqcup_{p \in P} P$. Let us denote the embedding of $P$ into the $p$-th copy of $P$ in the disjoint union by $P \ni x \ma p^{-1} \cdot x \in \bigsqcup_{p \in P} P$. Then we define an equivalence relation $\sim$ by identifying $p_1^{-1} \cdot x_1$ and $p_2^{-1} \cdot x_2$ in $\bigsqcup_{p \in P} P$ if there exists $p$ in $P$ with $p_1 \leq p$, $p_2 \leq p$ and $(pp_1^{-1})x_1 = (pp_2^{-1})x_2$. The set of equivalence classes $(\bigsqcup_{p \in P} P) / \sim$ with respect to $\sim$ carries the following canonical structure of a semigroup: Given $p_1^{-1} \cdot x_1$ and $p_2^{-1} \cdot x_2$ in $\bigsqcup_{p \in P} P$, take $y \in P$ with $x_1 \leq y$ and $p_2 \leq y$ and set
\bgl
\label{op:G}
  \eckl{p_1^{-1} \cdot x_1} \eckl{p_2^{-1} \cdot x_2} = \eckl{((yx_1^{-1})p_1)^{-1} \cdot ((yp_2^{-1})x_2)}.
\egl
Here $\eckl{\cdot}$ stands for equivalence class. One can check that the set $(\bigsqcup_{p \in P} P) / \sim$ together with the binary operation defined in \eqref{op:G} is indeed a semigroup. Let us denote it by $G$. The unit element in $G$ is given by $\eckl{e^{-1} \cdot e}$ where $e$ is the unit element of $P$. Moreover, by definition of the binary operation, we have
\bgloz
  \eckl{p^{-1} \cdot x} \eckl{x^{-1} \cdot p} = \eckl{p^{-1} \cdot p} = \eckl{e^{-1} \cdot e}
\egloz
(take $y = x$ in \eqref{op:G}). So we see that we have actually defined a group. Finally, the map $P \ni p \ma \eckl{e^{-1} \cdot p} \in G$ defines a semigroup homomorphism which is injective if $P$ is also left cancellative. By construction, the group $G$, together with this embedding of $P$, satisfies all the properties Theorem~\ref{OD}: Every element of $G$ is of the form $\eckl{p^{-1} \cdot x} = \eckl{p^{-1} \cdot e} \eckl{e^{-1} \cdot x} = (\eckl{e^{-1} \cdot p})^{-1} \eckl{e^{-1} \cdot x} \in P^{-1} P$. Here we are identifying $P$ with its image in $G$ under the embedding $P \ni p \ma \eckl{e^{-1} \cdot p} \in G$. Moreover, given a group $G'$ and a semigroup homomorphism $\varphi: P \to G'$, it is straightforward to check that the map $G \to G', \eckl{p^{-1} \cdot x} \ma \varphi(p)^{-1} \varphi (x)$ defines a group homomorphism which extends $\varphi$. Uniqueness of the extension follows from the equation $G = P^{-1} P$.

This is one way of constructing a model for the enveloping group. The main idea is to formally invert semigroup elements using an inductive limit procedure. Similar ideas frequently appear in the literature (compare for instance \cite{La}), and as we will see, this idea will also play a role later on in this paper.

\subsection{Reduced semigroup C*-algebras}

The main goal of this paper is to compute K-theory for reduced semigroup C*-algebras of left Ore semigroups whose constructible right ideals are independent (under a certain K-theoretic assumption on the enveloping group). In this paragraph, let us briefly recall the construction of reduced semigroup C*-algebras. The reader may consult \cite{Li2} for details.

Let $P$ be a left cancellative semigroup. Let $\ell^2(P)$ be the Hilbert space of square summable functions from $P$ to $\Cz$ and let $\menge{\varepsilon_x}{x \in P}$ be the canonical orthonormal basis of $\ell^2(P)$ given by $\varepsilon_x(y) = \delta_{x,y}$ ($\delta_{x,y} = 1$ if $x=y$ and $\delta_{x,y} = 0$ if $x \neq y$). The semigroup $P$ acts on $\ell^2(P)$ as follows: For every $p \in P$, the map $\varepsilon_x \ma \varepsilon_{px}$ extends to an isometry $V_p$ on $\ell^2(P)$ because our assumption that $P$ is left cancellative implies that $P \ni x \ma px \in P$ is injective. Now we simply set
\bsubdefin
$C^*_r(P) \defeq C^*(\menge{V_p}{p \in P}) \subseteq \cL(\ell^2(P))$.
\esubdefin
This is the reduced semigroup C*-algebra of $P$. In other words, the reduced semigroup C*-algebra is the C*-algebra generated by the left regular representation of the semigroup.

Now consider the family $\cJ$ of right ideals of $P$ from Definition~\ref{cJ}. For every right ideal $X \in \cJ$, we let $E_X$ be the orthogonal projection on $\ell^2(P)$ onto the subspace $\ell^2(X) \subseteq \ell^2(P)$. As observed in \cite{Li2}, \S~2, the projections $E_X$ lie in $C^*_r(P)$ for all $X \in \cJ$. Thus, the following
\bsubdefin
$D_r(P) \defeq C^*(\menge{E_X}{X \in \cJ}) \subseteq \cL(\ell^2(P))$

\textnormal{defines a sub-C*-algebra of $C^*_r(P)$. It is clear that $D_r(P)$ is a commutative C*-algebra, and that the multiplication on the generators is given by $E_X E_Y = E_{X \cap Y}$. Moreover, $D_r(P)$ is $\Ad(V_p)$-invariant for every $p \in P$. Therefore, the map $\tau: P \to \End(D_r(P)); p \ma \tau_p \defeq \Ad(V_p) \vert_{D_r(P)}$ defines a semigroup action of $P$ on $D_r(P)$.}
\esubdefin

\subsection{On reduced crossed products}

Let us collect a few observations about reduced crossed products. These results are included for the sake of completeness and also for ease of reference. They are certainly well known and we do not claim any originality here. In what follows we always assume that $G$ is a discrete group although most of what we say below has obvious analogues for general locally compact groups.

We denote by $\lambda:G\to \cU(\ell^2(G))$ the left regular representation of $G$ and by $M:c_0(G)\to \cL(\ell^2(G))$ the representation of $c_0(G)$ by multiplication operators on $\ell^2(G)$.

Recall that the reduced crossed product $A\rtimes_{\alpha,r}G$ of the C*-dynamical system $(A,G,\alpha)$ can be defined as the sub-C*-algebra of $M(A\otimes\cK_G)$, with $\cK_G:=\cK(\ell^2(G))$, generated by the set
$$\{\iota_A(a)\iota_G(g):a\in A, g\in G\},$$
where $\iota_G(g)=1\otimes \lambda_g$ and where $\iota_A:A\to M(A\otimes \cK_G)$ is defined by the composition
$$A\stackrel{\tilde{\alpha}}{\longrightarrow} \ell^{\infty}(G,A)\subseteq M(A\otimes c_0(G))\stackrel{\id_A\otimes M}{\longrightarrow} M(A\otimes \cK_G).$$
Here $\tilde\alpha$  sends $a\in A$ to the function $\eckl{g\mapsto \alpha_{g^{-1}}(a)}\in \ell^{\infty}(G,A)$. Every representation $\rho:A\to \cL(H)$ induces a homomorphism 
$$\Ind\rho: A\rtimes_rG\to \cL(H\otimes\ell^2(G))$$
by applying the representation $\rho\otimes\id_{\cK_G}: A\otimes \cK_G\to \cL(H\otimes\ell^2(G))$ to $A\rtimes_{\alpha,r}G\subseteq M(A\otimes \cK_G)$. It follows that $\Ind\rho$ is faithful if $\rho$ is faithful. One easily checks that
\begin{equation}\label{eq-Ind}
(\Ind\rho)(\iota_A(a))(\xi\otimes \varepsilon_g)=\rho(\alpha_{g^{-1}}(a))\xi\otimes\varepsilon_g\quad\text{and}\quad (\Ind\rho)(\iota_G(g))=1\otimes \lambda_g
\end{equation}
for all $a\in A$, $\xi\in H$ and $g\in G$, where $\{\varepsilon_x: x\in G\}$ denotes the standard orthonormal basis of $\ell^2(G)$. Thus, if $\rho:A\to \cL(H)$ is faithful, we recover the classical spatial definition of the reduced crossed product as a subalgebra of $\cL(H\otimes \ell^2(G))$.

Our first lemma is concerned with crossed products $D\rtimes_{\tau,r}G$ where $D$ is a closed, left-translation invariant sub-C*-algebra $D$ of $\ell^\infty(G)$ and $\tau:G\to \Aut(D)$ denotes the left-translation action. Let $M:D\to \cL(\ell^2(G))$ be the representation by multiplication operators. One easily checks that $(M,\lambda)$ is a covariant representation of $(D,G,\tau)$ on $\ell^2(G)$. It therefore induces a representation $M\rtimes\lambda: D\rtimes_\tau G\to \cL(\ell^2(G))$. 

\bsublemma
\label{DxG}
Let $(M,\lambda)$ be as above. Then $(M\otimes 1, \lambda\otimes 1)$ is unitarily equivalent to the regular representation 
$(\Ind M) \circ \iota_D, (\Ind M) \circ \iota_G)$ on $\ell^2(G\times G)$.
In particular,  $M \rtimes \lambda: D\rtimes_\tau G\to \cL(\ell^2(G))$ factors through a faithful representation of $D\rtimes_{\tau,r}G$.
\esublemma
\bproof
Consider the unitary operator $W: \ell^2(G\times G)\to \ell^2(G\times G); W(\varepsilon_x\otimes \varepsilon_y)=\varepsilon_{yx}\otimes \varepsilon_{x^{-1}}$; its adjoint is given by the formula $W^*(\varepsilon_x\otimes \varepsilon_y)=\varepsilon_{y^{-1}}\otimes \varepsilon_{xy}$. We then compute for $f \in \ell^\infty(G)$:
\begin{align*}
W ((\Ind M) \circ \iota_D)(f) W^*(\varepsilon_x\otimes \varepsilon_y)&=W ((\Ind M) \circ \iota_D)(f)(\varepsilon_{y^{-1}}\otimes \varepsilon_{xy})\\
&=W(f(x)(\varepsilon_{y^{-1}}\otimes \varepsilon_{xy})\\
&=f(x)(\varepsilon_x\otimes \varepsilon_y)=(M(f)\otimes 1)(\varepsilon_x\otimes \varepsilon_y)
\end{align*}
and
\begin{align*}
W(1\otimes\lambda_g)W^*(\varepsilon_x\otimes \varepsilon_y)&=
W(1\otimes \lambda_g)(\varepsilon_{y^{-1}}\otimes \varepsilon_{xy})\\
&=W(\varepsilon_{y^{-1}}\otimes \varepsilon_{gxy})=\varepsilon_{gx}\otimes \varepsilon_y\\
&=(\lambda_g \otimes 1)(\varepsilon_x\otimes \varepsilon_y).
\end{align*}
\eproof

Our second lemma is about functorial properties of reduced crossed products.
\bsublemma
\label{functor}
Suppose that $(A,G,\alpha)$ and $(B,H,\beta)$ are C*-dynamical systems, where $G$ and $H$ are discrete groups. Assume that $\varphi:A\to B$ is a homomorphism and that $j:G\to H$ is an injective homomorphism such that $\beta_{j(g)}(\varphi(a))=\varphi(\alpha_g(a))$ for all $a\in A$ and $g\in G$. Then there exists a unique homomorphism $\varphi\rtimes_rj:A\rtimes_{\alpha,r}G\to B\rtimes_{\beta,r}H$ such that 
$(\varphi\rtimes_rj)(\iota_A(a)\iota_G(g))=\iota_B(\varphi(a))\iota_H(j(g))$. If $\varphi$ is faithful, then so is $\varphi\rtimes_rj$.
\esublemma
\bproof We may assume without loss of generality that $G$ is a subgroup of $H$ and that $j:G\to H$ is the inclusion map. Restricting $\beta$ to $G$, we first observe that we have a homomorphism $\varphi \otimes \id_{\cK_G}: A\otimes\cK_G\to B\otimes \cK_G$ which we may extend to a homomorphism (again denoted by $\varphi \otimes \id_{\cK_G}$) $A\rtimes_{\alpha,r}G \to M(B\otimes\cK_G)$ such that 
$$(\varphi\otimes \id_{\cK_G})(\iota_A(a))=\iota_B(\varphi(a))\quad\text{and}\quad (\varphi\otimes \id_{\cK_G})(\iota_G(g))=\iota_G(g).$$
Thus $\varphi\otimes \id_{\cK_G}$ maps $A\rtimes_{\alpha,r}G$ into $B\rtimes_{\beta,r}G$ and $\varphi\otimes\id_{\cK_G}$ is faithful if $\varphi$ is faithful.

To see that $B\rtimes_{\beta,r}G$ imbeds into $B\rtimes_{\beta,r}H$, we first observe that $\ell^2(H)$ can be identified with $\bigoplus_{[h]\in G\backslash H}\ell^2(G)$. An explicit isomorphism is given by choosing a cross section $c:G\backslash H\to H$ which induces a bijection $G\times G\backslash H\to H; (g,[h])\mapsto gc([h])$ and hence an isomorphism $\ell^2(H)\cong \bigoplus_{[h]\in G\backslash H}\ell^2(G)$ by sending $\varepsilon_{gc[h]}$ to $\varepsilon_g$ in the summand at $[h]$ for all $g\in G, [h]\in G\backslash H$. Under this isomorphism, we get for $b\in B$ that $\iota_B^H(b)=
\bigoplus_{[h]\in G\backslash H}(\beta_{c[h]^{-1}}\otimes \id_{\cK_G})(\iota_B^G(b))\in \bigoplus_{[h]\in G\backslash H} M(B\otimes\cK_G)\subseteq M(B\otimes \cK_H)$ and $\iota_H(g)=\bigoplus_{[h]\in G\backslash H} \iota_G(g)$ for all $g\in G$ (where the superscript $H$ indicates that $\iota_B^H(b)$ belongs to the crossed product $B\rtimes_{\beta,r}H$). Thus we see that the subalgebra of $B\rtimes_{\beta,r}H$ generated by $\{\iota_B^H(b)\iota_H(g): b\in B, g\in G\}$ equals
$\bigoplus_{[h]\in G\backslash H} (\beta_{c[h]^{-1}}\otimes\id_{\cK_G})(B\rtimes_{\beta,r}G)$
which is isomorphic to $B\rtimes_{\beta,r}G$ via an isomorphism sending $\iota_B^G(b)\iota_G(g)$ to $\iota_B^H(b)\iota_H(g)$ for all $b\in B$, $g\in G$. Combining this with the first part gives the lemma.
\eproof

For the proof of the following lemma we refer to \cite{Br-Oz}, Chapter~4, Proposition~1.9.
\bsublemma
\label{fce}
Let $(A,G,\alpha)$ be a C*-dynamical system with $G$ discrete. Then there exists a unique faithful conditional expectation $E:A\rtimes_{\alpha,r}G\to A$ such that $E(\iota_A(a)\iota_G(g))=\delta_{g,e}a$, where $\delta_{g,e}=1$ if $g$ is equal to the unit $e$ of $G$ and $\delta_{g,e}=0$ if $g\neq e$.
\esublemma

\section{The strategy}

Let $P$ be a left Ore semigroup whose constructible right ideals are independent. Let $G$ be the enveloping group of $P$ (see Theorem~\ref{OD}). Using Theorem~\ref{OD}, we will always view $P$ as a subsemigroup of $G$.

Our goal is to compute K-theory for the reduced semigroup C*-algebra of $P$ under a K-theoretic assumption on $G$ which we will make precise later on. Let us now present our strategy:

First, we make use of the assumption that $P$ is a left Ore semigroup to reduce our K-theoretic problem to the problem of computing K-theory for a reduced crossed product by the enveloping group $G$ of $P$. The main idea has already appeared in the previous section, namely to use inductive limit procedures to pass from $P$ to $G$.

The main step is to compare the reduced crossed product we are interested in with another, but much simpler reduced crossed product. The simpler one is given by an action of $G$ on a discrete space (simply a set). This step makes use of our K-theoretic assumption on $G$. It allows us to apply the machinery of Baum-Connes which will reduce the K-theoretic comparison of the reduced crossed products to the case of finite subgroups. Here our assumption that the constructible right ideals of $P$ are independent enters the game, as we will see.

The last step is to compute K-theory for reduced crossed products associated with an action of our group $G$ on a discrete space. This amounts to applying imprimitivity theorems.

\section{Dilations of reduced semigroup C*-algebras}
\label{dilations}

For what we are going to do in this section, it is enough to assume that our semigroup $P$ satisfies the left Ore condition. We would like to describe the reduced semigroup C*-algebra $C^*_r(P)$ as a reduced crossed product by the enveloping group $G$, at least up to Morita equivalence. Following ideas of \cite{La}, we first of all construct a $G$-C*-algebra which gives rise to the reduced crossed product.

Similarly as in Paragraph~\ref{Ore}, we consider the following inductive system of C*-algebras indexed by elements of $P$ ordered by \an{$\leq$}:
\begin{itemize}
\item the $p$-th semigroup is given by $D_r(P)$ for every $p \in P$
\item given $p, q \in P$ with $p \leq q$, the structure map from the $p$-th to the $q$-th C*-algebra is given by $\tau_{q p^{-1}} = \Ad(V_{q p^{-1}}): D_r(P) \to D_r(P)$.
\end{itemize}
Let $\DiP$ be the inductive limit of this system, and denote by $\iota_p: \DP \to \DiP$ the inclusion of the $p$-th C*-algebra into the inductive limit. As explained in \cite{La}, there is a $G$-action $\taui$ on $\DiP$ which dilates the $P$-action $\tau$ on $\DP$. To describe $\taui$, it suffices to define $\taui_p$ for every $p \in P \subseteq G$ because the semigroup homomorphism $P \ni p \ma \taui_p \in \Aut(\DiP)$ extends uniquely to $G$ by Theorem~\ref{OD}. Now $\taui_p$ is given as follows: For $q \in P$ and $d \in \DP$, let $r$ be an element in $P$ such that $p \leq r$ and $q \leq r$. Then we set
\bgloz
  \taui_p (\iota_q(d)) \defeq \iota_{rp^{-1}}(\tau_{rq^{-1}}(d)).
\egloz
One can check that this formula gives rise to the desired automorphism $\taui_p$ of $\DiP$ and that these automorphisms give rise to the semigroup homomorphism $P \to \Aut(\DiP)$, $p \ma \taui_p$. Moreover, one can also verify that the automorphisms we have constructed coincide with the ones constructed in \S~2 of \cite{La}.

In the following, we construct a covariant representation for the C*-dynamical system $(\DiP, G, \taui)$. First, we obtain a canonical faithful representation of $\DiP$ on $\ell^2(G)$ as follows: Using the inductive limit structure of $\DiP$, it suffices to construct a family of faithful representations $\gekl{\pi_p}_{p \in P}$ of $\DP$ on $\ell^2(G)$ which are compatible with the structure maps. As $\DP$ acts on $\ell^2(P)$ by construction, we can conjugate the identity representation of $\DP$ by the canonical isometric embedding $\ell^2(P) \into \ell^2(G)$ to obtain a faithful representation $\pi$ of $\DP$ on $\ell^2(G)$. Then define for every $p \in P$ the representation $\pi_p \defeq \Ad(\lambda_p^*) \circ \pi$. Here for every $g \in G$, we denote by $\lambda_g$ the unitary on $\ell^2(G)$ given by $\lambda_g(\varepsilon_x) = \varepsilon_{gx}$ for the canonical orthonormal basis $\menge{\varepsilon_x}{x \in G}$ of $\ell^2(G)$. In other words, $\lambda_g$ is the image of $g \in G$ under the left regular representation $\lambda$ of $G$. These representations $\pi_p$ are faithful by construction. For a subset $Y$ of $G$, let $E_Y \in \cL(\ell^2(G))$ be the orthogonal projection onto the subspace $\ell^2(Y)$ of $\ell^2(G)$. It is then immediate that for every $X \in \cJ$, we have 
\bgl
\label{piEX}
  \pi_p(E_X) = E_{p^{-1} \cdot X}.
\egl
Note that $p^{-1} \cdot X$ is the subset $\menge{p^{-1} x}{x \in X}$ of $G$; it should not be confused with $p^{-1} X = \menge{q \in P}{pq \in X}$. From \eqref{piEX}, it follows that the representations $\pi_p$ are compatible with the structure maps, in the sense that for every $p, q \in P$ with $p \leq q$, we have
\bgloz
  \pi_q \circ \Ad(V_{qp^{-1}}) = \pi_p.
\egloz
Therefore, the faithful representations $\gekl{\pi_p}_{p \in P}$ give rise to a faithful representation $\pii$ of $\DiP$ on $\ell^2(G)$. This representation is determined by $\pii(\iota_q(E_X)) = E_{q^{-1} \cdot X}$.

We claim that this representation $\pii$, together with the left regular representation $\lambda$ of $G$, is a covariant representation of $(\DiP, G, \taui)$. To show this, take $p, q, r \in P$ with $p \leq r$, $q \leq r$, $X \in \cJ$ and compute
\bglnoz
  \lambda_p (\pii(\iota_q(E_X))) \lambda_p^* &=& \lambda_p E_{q^{-1} \cdot X} \lambda_p^* = E_{p \cdot q^{-1} \cdot X} \\
  &=& E_{p \cdot r^{-1} \cdot r \cdot q^{-1} \cdot X} = E_{(r p^{-1})^{-1} \cdot (r q^{-1}) \cdot X} \\
  &=& \pii(\iota_{rp^{-1}}(\tau_{rq^{-1}}(E_X))) = \pii(\taui_p(\iota_q(E_X))).
\eglnoz
So far, we have constructed a covariant representation $(\pii,\lambda)$ of the C*-dynamical system $(\DiP,G,\taui)$. Next we claim:
\blemma
\label{redcp-rep}
The covariant representation $(\pii,\lambda)$ gives rise to a faithful representation $\pii \rtimes_r \lambda$ of the reduced crossed product $\DiP \rtimes_{\taui,r} G$ on $\ell^2(G)$. This representation is determined by $(\pii \rtimes_r \lambda)(d U_g) = \pii(d) \lambda_g$ for all $d \in \DiP$ and $g \in G$. Here $U_g$ are the canonical unitaries in the multiplier algebra of $\DiP \rtimes_{\taui,r} G$ implementing $\taui$.
\elemma
\bproof
Apply Lemma~\ref{DxG} to $D = \pii(\DiP)$.
\eproof
Using this representation $\pii \rtimes_r \lambda$, we will always think of $\DiP \rtimes_{\taui,r} G$ as a concrete C*-algebra acting on $\ell^2(G)$.

Now consider the orthogonal projection $E_P \in \cL(\ell^2(G))$ onto the subspace $\ell^2(P) \subseteq \ell^2(G)$. This projection lies in $\DiP \rtimes_{\taui,r} G$.
\blemma
\label{P-EGE}
The projection $E_P$ is a full projection in $\DiP \rtimes_{\taui,r} G$, and the corner $E_P \rukl{\DiP \rtimes_{\taui,r} G} E_P$ can be identified with $C^*_r(P)$ via
\bgloz
  C^*_r(P) \ni V_p \ma E_P U_p E_P \in E_P \rukl{\DiP \rtimes_{\taui,r} G} E_P.
\egloz
\elemma
\bproof
The C*-algebra $\DiP$ (or rather $\pii(\DiP)$) is generated by the projections $\menge{E_{q^{-1} \cdot X}}{q \in P, X \in \cJ}$. Thus the net $(E_{q^{-1} \cdot P})_{q \in P}$ is an approximate unit of $\DiP$, hence of $\DiP \rtimes_{\taui,r} G$. As
\bgloz
  \E{q}{P} = U_q^* E_P U_q \in \rukl{\DiP \rtimes_{\taui,r} G} E_P \rukl{\DiP \rtimes_{\taui,r} G},
\egloz
our first claim follows.

Now let us prove that the assignment $V_p \ma E_P U_p E_P$ extends to an isomorphism
\bgloz
  C^*_r(P) \to E_P \rukl{\DiP \rtimes_{\taui,r} G} E_P.
\egloz
The assignment $V_p \ma E_P U_p E_P$ first of all extends to a homomorphism $C^*_r(P) \to E_P \rukl{\DiP \rtimes_{\taui,r} G} E_P$ because the operator $E_P U_p E_P$, viewed as an operator on $\ell^2(P) \subseteq \ell^2(G)$, is really nothing else but the isometry $V_p$ itself (note that $U_p$ is just $\lambda_p$ since we view $\DiP \rtimes_{\taui,r} G$ as a concrete C*-algebra acting on $\ell^2(G)$ via $\pii \rtimes_r \lambda$). This observation implies that the resulting homomorphism $C^*_r(P) \to E_P \rukl{\DiP \rtimes_{\taui,r} G} E_P$ must be injective. To show surjectivity, it is enough to prove that for all $p, q_1, q_2 \in P$ and $X \in \cJ$, the element $E_P \E{q_1}{X} U_{q_2^{-1}p} E_P \in E_P \rukl{\DiP \rtimes_{\taui,r} G} E_P$ lies in the image. But
\bglnoz
  E_P \E{q_1}{X} U_{q_2^{-1}p} E_P &=& \rukl{E_P \E{q_1}{X} E_P} \rukl{E_P U_{q_2}^* E_P} \rukl{E_P U_p E_P} \\
  &=& \rukl{E_P E_{P \cap (q_1^{-1} \cdot X)} E_P} \rukl{E_P U_{q_2} E_P}^* \rukl{E_P U_p E_P} \\
  &=& \rukl{E_P \E{q_1}{X} E_P} \rukl{E_P U_{q_2} E_P}^* \rukl{E_P U_p E_P}
\eglnoz
is the image of $E_{q_1^{-1} X} V_{q_2}^* V_p$.
\eproof
\bcor
\label{CP-DG}
The embedding $\iota: C^*_r(P) \to \DiP \rtimes_{\taui,r} G$ determined by $\iota(V_p) =  E_P U_p E_P$ induces a KK-equivalence in $KK(C^*_r(P), \DiP \rtimes_{\taui,r} G)$.
\ecor

\section{From concrete to abstract}
\label{concrete-abstract}

Corollary~\ref{CP-DG} tells us that if we are interested in the K-theory of $C^*_r(P)$, we can equally well study the reduced crossed product $\DiP \rtimes_{\taui,r} G$. The situation is as follows:

\begin{itemize}
\item[(i)] $\DiP$ is a commutative C*-algebra generated by the projections
\bgloz
  \menge{\E{q}{X}}{q \in P, \: \emptyset \neq X \in \cJ}.
\egloz
As $P$ is countable, this family of projections is countable as well ($\cJ$ is countable as $P$ is). Moreover,
\bgloz
  \menge{\E{q}{X}}{q \in P, X \in \cJ} = \menge{\E{q}{X}}{q \in P, \emptyset \neq X \in \cJ} \cup \gekl{0}
\egloz
is multiplicatively closed because given $q_1, q_2 \in P$ and $X_1, X_2 \in \cJ$, we can choose $q \in P $ with $q_1 \leq q$ and $q_2 \leq q$, and then
\bgloz
  (q_1^{-1} \cdot X_1) \cap (q_2^{-1} \cdot X_2) = q^{-1} \cdot \underbrace{\rukl{(qq_1^{-1} \cdot X_1) \cap (qq_2^{-1} \cdot X_2)}}_{\in \cJ}
\egloz
so that
\bgloz
  \E{q_1}{X_1} \E{q_2}{X_2} = E_{q^{-1} \cdot \rukl{(qq_1^{-1} \cdot X_1) \cap (qq_2^{-1} \cdot X_2)}}
\egloz
lies in $\menge{\E{q}{X}}{q \in P, X \in \cJ}$.

\item[(ii)] Assume that the constructible right ideals of $P$ are independent. Then we can prove the following:

For all projections $E$, $E_1, \dotsc, E_n$ in $\menge{\E{q}{X}}{q \in P, X \in \cJ}$, the strict inequalities $E_1, \dotsc, E_n \lneq E$ imply $\bigvee_{j=1}^n E_j \lneq E$. Here $\bigvee_{j=1}^n E_j$ is the smallest projection in $\DiP$ which is bigger than or equal to $E_1, \dotsc, E_n$.

Here is the proof: Let $E_j = \E{q_j}{X_j}$, $j = 1, \dotsc, n$, and $E = \E{q}{X}$ with $q, q_1, \dotsc, q_n \in P$ and $X, X_1, \dotsc, X_n \in \cJ$. It follows that $\bigvee_{j=1}^n \E{q_j}{X_j} = E_{\bigcup_{j=1}^n q_j^{-1} \cdot X_j}$. We claim that $\E{q_j}{X_j} \lneq \E{q}{X}$ for all $1 \leq j \leq n$ implies $\bigvee_{j=1}^n \E{q_j}{X_j} \lneq \E{q}{X}$. As for $Y_1, Y_2 \subseteq G$, the inequality $E_{Y_1} \leq E_{Y_2}$ is equivalent to $Y_1 \subseteq Y_2$, we have to show that $q_j^{-1} \cdot X_j \subsetneq q^{-1} \cdot X$ for all $1 \leq j \leq n$ implies $\bigcup_{j=1}^n q_j^{-1} \cdot X_j \subsetneq q^{-1} \cdot X$. But this follows from our assumption that the constructible right ideals of $P$ are independent: Choose $r \in P$ such that $q \leq r$ and $q_j \leq r$ for all $1 \leq j \leq n$. Then $q_j^{-1} \cdot X_j \subsetneq q^{-1} \cdot X$ for all $1 \leq j \leq n$ implies that $(r q_j^{-1}) \cdot X_j \subsetneq (r q^{-1}) \cdot X$ for all $1 \leq j \leq n$. But $(r q^{-1}) \cdot X = (r q^{-1}) X$ and $(r q_j^{-1}) \cdot X_j = (r q_j^{-1}) X_j$ lie in $\cJ$ for all $1 \leq j \leq n$. Thus the independence condition tells us that
\bgloz
  r \rukl{\bigcup_{j=1}^n q_j^{-1} \cdot X_j} = \bigcup_{j=1}^n (rq_j^{-1}) X_j \subsetneq (rq^{-1}) X = r \rukl{q^{-1} \cdot X}.
\egloz
Since left multiplication by $r$ is injective, we deduce $\bigcup_{j=1}^n q_j^{-1} \cdot X_j \subsetneq q^{-1} \cdot X$, as claimed.

\item[(iii)] The $G$-action $\taui$ on $\DiP$ leaves the projections
\bgloz
  \menge{\E{q}{X}}{q \in P, \emptyset \neq X \in \cJ}
\egloz
invariant.
\end{itemize}

This is the situation we are interested in. In the following section, we look at it from an abstract point of view.

\section{The general K-theoretic result}
\label{generalK}

We first formulate our assumptions:
\begin{itemize}
\item[(I)]
$D$ is a commutative C*-algebra generated by a countable family of pairwise distinct (commuting) non-zero projections $\gekl{e_i}_{i \in I}$. Moreover, $\gekl{e_i}_{i \in I} \cup \gekl{0}$ is multiplicatively closed (i.e. for all $e_i$, $e_j$ in $\gekl{e_i}_{i \in I}$, either $e_i e_j = 0$ or there exists $k \in I$ such that $e_i e_j = e_k$).

\item[(II)]
The family $\gekl{e_i}_{i \in I}$ is independent, i.e. given $e \in \gekl{e_i}_{i \in I}$ and finitely many $e_1, \dotsc, e_n \in \gekl{e_i}_{i \in I}$ with $e_1, \dotsc, e_n \lneq e$, we always have $\bigvee_{i=1}^n e_i \lneq e$, i.e. $e - \bigvee_{i=1}^n e_i$ is a non-zero projection. Here $\bigvee_{i=1}^n e_i$ is the smallest projection in $D$ which is bigger than (or equal to) all the $e_i$, $1 \leq i \leq n$. Note that since $D$ is commutative, $\bigvee_{i=1}^n e_i = \sum_{\emptyset \neq J \subseteq \gekl{1, \dotsc, n}} (-1)^{\abs{J}-1} \prod_{j \in J} e_j$. 

\item[(III)]
$G$ is a discrete countable group and $\tau$ is an action of $G$ on $D$ which leaves $\gekl{e_i}_{i \in I}$ invariant. This means that there is an action of $G$ on the index set $I$ such that $\tau_g(e_i) = e_{g \cdot i}$.
\end{itemize}

Assume that we have a C*-dynamical system $(D,G,\tau)$ satisfying (I), (II) and (III). In this situation, the homomorphisms $\phi_i: \Cz \to D; 1 \ma e_i$ (for $i \in I$) give rise to a KK-element in $KK(\bigoplus_{i \in I} \Cz,D) \cong \prod_{i \in I} KK(\Cz,D)$ which can be viewed as an element in equivariant KK-theory. This means that with respect to the $G$-action $\sigma$ on $\bigoplus_{i \in I} \Cz$ given by shifting the index set $I$ and the $G$-action $\tau$ on $D$, the $\phi_i$ yield in a canonical way an element $\xf \in KK^G(\bigoplus_{i \in I} \Cz,D)$. This KK-element will be described in detail in \S~\ref{proofmainthm-beginning}. Here is our main result:

\btheo
\label{mainthm}
Let us assume that we are in the situation described above. Then for every finite subgroup $H$ of $G$, the element $j^H(\res_H^G(\xf)) \in KK((\bigoplus_{i \in I} \Cz) \rtimes_{\sigma} H, D \rtimes_{\tau} H)$ is a KK-equivalence. Here $\res_H^G$ is the canonical restriction map $KK^G \to KK^H$ and $j^H$ is the descent $KK^H(\bigoplus_{i \in I} \Cz,D) \to KK((\bigoplus_{i \in I} \Cz) \rtimes_{\sigma} H, D \rtimes_{\tau} H)$.
\etheo

The proof of this theorem is the content of \S~\ref{proofmainthm-beginning} to \S~\ref{proofmainthm-end}.

Just a remark on notation: From now on, we write $c_0(I)$ for $\bigoplus_{i \in I} \Cz$ and $c_0(I,D)$ for $\bigoplus_{i \in I} D$. Moreover, given a Hilbert module $Z$, we write $\ell^2(I,Z)$ for $\bigoplus_{i \in I} Z$ (where the direct sum is taken in the sense of Hilbert modules).

\subsection{Description of the KK-element}
\label{proofmainthm-beginning}

Our goal is to describe the element $\xf \in KK^G(c_0(I),D)$. First of all, the element in $KK(c_0(I),D)$ given by the homomorphisms $\phi_i: \Cz \to D; 1 \ma e_i$ ($i \in I$) can be represented by the Kasparov module $(\ell^2(I,D), \phi, 0)$. The left action of $c_0(I)$ on the Hilbert D-module $\ell^2(I,D)$ is given by
\bgloz
  \phi \defeq \bigoplus \phi_i: c_0(I) \to c_0(I,D) \subseteq \cK(\ell^2(I,D)) \subseteq \cL(\ell^2(I,D)).
\egloz
Here $c_0(I,D)$ acts on $\ell^2(I,D)$ by diagonal operators. Let us introduce the notation that $\1z_i \in c_0(I)$ is the element whose $i$-th component is $1$ and whose other components are $0$, and we denote by $\1z_j \otimes d \in \ell^2(I,D)$ the element whose $j$-th component is $d \in D$ and whose remaining components vanish. Then
\bgloz
  \phi(\1z_i)(\1z_j \otimes d) = (\1z_i \otimes \phi_i(1))(\1z_j \otimes d) = (\1z_i \otimes e_i)(\1z_j \otimes d) = \delta_{i,j} (\1z_i \otimes e_i d).
\egloz
Since $\img(\phi)$ is contained in the compact operators on $\ell^2(I,D)$, the operator in our Kasparov module may be chosen to be $0$.

Now we want to interpret this Kasparov module as an element in $KK^G(c_0(I),D)$. So we introduce a $G$-action on the Hilbert module $\ell^2(I,D)$ which is compatible with the action $\tau$ of $G$ on $D$ so that $\phi$ becomes $G$-equivariant.

We let $\sigma$ be the $G$-action on $c_0(I)$ determined by $\sigma_g(\1z_i) = \1z_{g \cdot i}$. The $G$-action $\sigma \otimes \tau$ on the Hilbert module $\ell^2(I,D)$ is given by
\bgloz
  (\sigma \otimes \tau)_g (\1z_i \otimes d) = \1z_{g \cdot i} \otimes \tau_g(d).
\egloz
Recall that $\1z_i \otimes d \in \ell^2(I,D)$ is the element whose $i$-th component is $d \in D$ and whose remaining components vanish. It can be checked immediately that this $G$-action $\sigma \otimes \tau$ is compatible with the Hilbert $D$-module structure on $\ell^2(I,D)$, in the sense that $\spkl{(\sigma \otimes \tau)_g (\xi), (\sigma \otimes \tau)_g (\eta)}_D = \tau_g(\spkl{\xi,\eta}_D)$ and $(\sigma \otimes \tau)_g (\xi \cdot d) = (\sigma \otimes \tau)_g (\xi) \cdot \tau_g(d)$ for all $g \in G$; $\xi$, $\eta$ in $\ell^2(I,D)$ and $d \in D$. Conjugation yields a $G$-action $\Ad(\sigma \otimes \tau)$ of $G$ on $\cL(\ell^2(I,D))$ given by $G \ni g \ma \Ad((\sigma \otimes \tau)_g) \in \Aut(\cL(\ell^2(I,D)))$. To check that $\phi$ is $G$-equivariant with respect to the $G$-action $\sigma$ on $c_0(I)$ and $\Ad(\sigma \otimes \tau)$, it suffices to consider elements $\1z_i \in c_0(I)$ and $\1z_j \otimes d \in \ell^2(I,D)$. We compute
\bglnoz
  && \rukl{\phi (\sigma_g(\1z_i))}(\1z_j \otimes d) = \rukl{\phi(\1z_{g \cdot i})}(\1z_j \otimes d) = (\1z_{g \cdot i} \otimes e_{g \cdot i})(\1z_j \otimes d) \\
  &=& \delta_{g \cdot i, j} \1z_{g \cdot i} \otimes (e_{g \cdot i} d) = \delta_{i, g^{-1} \cdot j} \1z_{g \cdot i} \otimes (e_{g \cdot i} d) 
  = (\sigma \otimes \tau)_g \rukl{\delta_{i, g^{-1} \cdot j} \1z_i \otimes (e_i \tau_{g^{-1}}(d))} \\
  &=& \Ad(\sigma \otimes \tau)_g \rukl{\phi(\1z_i)} (\1z_j \otimes d).
\eglnoz
This shows that the Kasparov module $(\ell^2(I,D), \phi, 0)$ together with the $G$-action $\sigma \otimes \tau$ really gives rise to an element $\xf \in KK^G(c_0(I),D)$.

Let us summarize our construction in the following
\bsubdefin
\label{def-xf}
Let $\xf$ be the element in $KK^G(c_0(I),D)$ (where $G$ acts on $c_0(I)$ and $D$ via $\sigma$ and $\tau$) which is represented by the Kasparov $G$-module for $(c_0(I),D)$ consisting of
\begin{itemize}
\item the Hilbert $D$-module $\ell^2(I,D)$ with $G$-action $\sigma \otimes \tau$ given by 
\bgloz
  (\sigma \otimes \tau)_g (\1z_i \otimes d) = \1z_{g \cdot i} \otimes \tau_g(d)
\egloz
\item the equivariant homomorphism 
\bgloz
  \phi: c_0(I) \to \cK(\ell^2(I,D)) \subseteq \cL(\ell^2(I,D))
\egloz
determined by $\rukl{\phi(\1z_i)}(\1z_j \otimes d) = \delta_{i,j} \1z_i \otimes e_i d$
\item the operator $0 \in \cL(\ell^2(I,D))$.
\end{itemize}
\esubdefin

\subsection{Descent of the restriction}
\label{jresx}

Let $H \subseteq G$ be a subgroup. Our goal is to describe the element $j_r^H(\res_H^G(\xf)) \in KK(c_0(I) \rtimes_{\sigma,r} H, D \rtimes_{\tau,r} H)$ given by the descent of the restriction of $\xf$ to $H$.

\bsubprop
For every subgroup $H$ of $G$, the KK-element $j_r^H(\res_H^G(\xf))$ in $KK(c_0(I) \rtimes_{\sigma,r} H, D \rtimes_{\tau,r} H)$ is represented by the Kasparov $(c_0(I) \rtimes_{\sigma,r} H, D \rtimes_{\tau,r} H)$-module consisting of
\begin{itemize}
\item the Hilbert $D \rtimes_{\tau,r} H$-module $\ell^2(I,(D \rtimes_{\tau,r} H))$
\item the homomorphism $\phi \rtimes_r H: c_0(I) \rtimes_{\sigma,r} H \to \cL(\ell^2(I,(D \rtimes_{\tau,r} H)))$ given by $(\phi \rtimes_r H)(\1z_i U_h) = (\1z_i \otimes e_i) \circ (\sigma_h \otimes U_h)$
\item the operator $0 \in \cL(\ell^2(I,(D \rtimes_{\tau,r} H)))$.
\end{itemize}
Here $U_h$ are the canonical unitaries in the multiplier algebra of $c_0(I) \rtimes_{\sigma,r} H$ which implement $\sigma$.
\esubprop
\bproof
First of all, to obtain a Kasparov $H$-module for $(c_0(I),D)$ with respect to the restricted actions $\sigma \vert_H$ and $\tau \vert_H$ (we will again denote these actions by $\sigma$ and $\tau$ in the sequel) which represents $\res_H^G(\xf)$, we can just take the Kasparov $G$-module from Definition~\ref{def-xf} and restrict the $G$-action $\sigma \otimes \tau$ to $H$.

We now describe the element $j_r^H(\res_H^G(\xf)) \in KK(c_0(I) \rtimes_{\sigma,r} H, D \rtimes_{\tau,r} H)$ following \cite{Kas}, \S~3.7. The construction for full crossed products is also described in \cite{Bla}, Chapter~VIII, \S~20.6, and it is very similar to the one for reduced crossed products. Of course, in the case of finite subgroups (which is in view of Theorem~\ref{mainthm} the most interesting case), it does not matter at all whether we take full or reduced crossed products.

By definition, $j_r^H(\res_H^G(\xf))$ is represented by the Kasparov module
\bgloz
  (\ell^2(I,D) \rtimes_{\tau,r} H, \psi, 0).
\egloz
Let us start with the Hilbert $D \rtimes_{\tau,r} H$-module $\ell^2(I,D) \rtimes_{\tau,r} H$. It is the completion of the pre-Hilbert $C_c(H,D)$-module whose underlying vector space $C_c(H,\ell^2(I,D))$ is given by all functions from $H$ to $\ell^2(I,D)$ with finite support (we are in the discrete case). Given such a function $\xi \in C_c(H,\ell^2(I,D))$ and an element $b \in C_c(H,D)$, the right action of $C_c(H,D)$ on $C_c(H,\ell^2(I,D))$ is given by
\bgl
\label{rightaction}
  (\xi \bullet b)(h) = \sum_{\ti{h} \in H} \xi(\ti{h}) \tau_{\ti{h}}(b(\ti{h}^{-1} h)).
\egl
Given two functions $\xi, \eta \in C_c(H,\ell^2(I,D))$, the $C_c(H,D)$-valued inner product on $C_c(H,\ell^2(I,D))$ is given by
\bgl
\label{innerproduct}
  \spkl{\xi,\eta}_{C_c(H,D)}(h) = \sum_{\ti{h} \in H} \tau_{\ti{h}^{-1}} (\spkl{\xi(\ti{h}),\eta(\ti{h}h)}_D).
\egl
Consider $D \rtimes_{\tau,r} H$ as a Hilbert module over itself and form the direct sum $\ell^2(I,(D \rtimes_{\tau,r} H))$. We claim that the map
\bgloz
  \Theta: C_c(H,\ell^2(I,D)) \to \ell^2(I,(D \rtimes_{\tau,r} H)); \: \xi \ma (\eckl{h \ma (\xi(h))_i})_i
\egloz
extends to an isomorphism $\Theta: \ell^2(I,D) \rtimes_{\tau,r} H \overset{\cong}{\lori} \ell^2(I,(D \rtimes_{\tau,r} H))$ of Hilbert $D \rtimes_{\tau,r} H$-modules. Here we view functions from $H$ to $D$ with finite support as elements of $D \rtimes_{\tau,r} H$, and a function $f: H \to D$ is often denoted by $\eckl{h \ma f(h)}$.

As $\Theta$ obviously has dense image in $\ell^2(I,(D \rtimes_{\tau,r} H))$, it suffices to check that $\Theta$ preserves the right $D \rtimes_{\tau,r} H$-actions as well as the inner products. It certainly suffices to check this for elements in $C_c(H,\ell^2(I,D))$ of the form $(\1z_i \otimes d) U_h = \eckl{\ti{h} \ma \delta_{\ti{h},h} \1z_i \otimes d}$. Such an element corresponds to $\1z_i \otimes (d U_h) \in \ell^2(I,(D \rtimes_{\tau,r} H))$ under $\Theta$. Here and in the sequel, $U_h$ is the characteristic function of $h \in H$. For the right $D \rtimes_{\tau,r} H$-actions, it certainly suffices to look at elements in $C_c(H,D) \subseteq D \rtimes_{\tau,r} H$ of the form $b = d_b U_{h_b}$. For $\xi = (\1z_i \otimes d_\xi) U_{h_\xi} \in C_c(H,\ell^2(I,D))$ and $b = d_b U_{h_b}$, we have by \eqref{rightaction}
\bglnoz
  (\xi \bullet b)(h) &=& \sum_{\tih \in H} \delta_{\tih,h_\xi} (\1z_i \otimes d_\xi) \tau_{\tih}(\delta_{\tih^{-1}h,h_b} d_b) 
  = \delta_{h_\xi^{-1}h,h_b} (\1z_i \otimes d_\xi)(\tau_{h_\xi}(d_b)) \\
  &=& \delta_{h_\xi^{-1}h,h_b} \1z_i \otimes (d_\xi \tau_{h_\xi}(d_b)) = (\1z_i \otimes d_\xi \tau_{h_\xi}(d_b)) U_{h_\xi h_b} (h)
\eglnoz
so that
\bglnoz
  \Theta(\xi \bullet b) &=& \1z_i \otimes (d_\xi \tau_{h_\xi}(d_b) U_{h_\xi h_b}) = \1z_i \otimes (d_\xi U_{h_\xi})(d_b U_{h_b}) 
  = (\1z_i \otimes (d_\xi U_{h_\xi}))(d_b U_{h_b}) \\
  &=& (\Theta(\xi)) \cdot b
\eglnoz
where in the last line, we let $b$ act on $\Theta(\xi) \in \ell^2(I,(D \rtimes_\tau H))$ using the right $D \rtimes_\tau H$-module structure of $\ell^2(I,(D \rtimes_\tau H))$.

Moreover, for $\xi = (\1z_i \otimes d_\xi) U_{h_\xi}$ and $\eta = (\1z_j \otimes d_\eta) U_{h_\eta}$ in $C_c(H,\ell^2(I,D))$, we have by \eqref{innerproduct}
\bglnoz
  \spkl{\xi,\eta}_{C_c(H,D)}(h) 
  &=& \sum_{\tih \in H} \tau_{\tih^{-1}}(\spkl{\delta_{\tih,h_\xi} \1z_i \otimes d_\xi, \delta_{\tih h,h_\eta} \1z_j \otimes d_\eta}_D) \\
  &=& \delta_{h_\xi h,h_\eta} \tau_{h_\xi^{-1}} (\delta_{i,j} d_\xi^* d_\eta) = \delta_{i,j} \delta_{h,h_\xi^{-1} h_\eta} \tau_{h_\xi^{-1}} (d_\xi^* d_\eta) \\
  &=& \delta_{i,j} (d_\xi U_{h_\xi})^* (d_\eta U_{h_\eta}) (h) \\
  &=& \spkl{\1z_i \otimes (d_\xi U_{h_\xi}), \1z_j \otimes (d_\eta U_{h_\eta})}_{D \rtimes_{\tau,r} H} (h) \\
  &=& \spkl{\Theta(\xi),\Theta(\eta)}_{D \rtimes_{\tau,r} H} (h).
\eglnoz
This proves our claim that $\Theta$ extends to an isomorphism of Hilbert $D \rtimes_{\tau,r} H$-modules.

Finally, it remains to describe $\psi$, i.e. to describe the left $c_0(I) \rtimes_{\sigma,r} H$-action on the Hilbert module. Let $a$ be an element in $C_c(H,c_0(I)) \subseteq c_0(I) \rtimes_{\sigma,r} H$. Then for $\xi \in C_c(H,\ell^2(I,D))$, $\psi(a) \xi$ is given by
\bgloz
  (\psi(a) \xi)(h) = \sum_{\tih \in H} \phi(a(\tih)) \rukl{(\sigma \otimes \tau)_{\tih} \xi(\tih^{-1} h)}
\egloz
(see \cite{Kas}, \S~3.7 or \cite{Bla}, Chapter~VIII, \S~20.6). To explicitly compute the action, we again take $\xi = (\1z_j \otimes d_\xi) U_{h_\xi}$ and $a = \1z_i U_{h_a}$. Then
\bglnoz
  (\psi(a) \xi)(h) &=& \sum_{\tih \in H} \delta_{\tih,h_a} \phi(\1z_i) \rukl{(\sigma \otimes \tau)_{\tih} \delta_{\tih^{-1}h,h_\xi} \1z_j \otimes d_\xi} \\
  &=& \delta_{h_a^{-1}h,h_\xi} \phi_i(\1z_i) (\1z_{h_a \cdot j} \otimes \tau_{h_a}(d_\xi)) \\
  &=& \delta_{h_a^{-1}h,h_\xi} \delta_{i,h_a \cdot j} (\1z_i \otimes e_i \tau_{h_a}(d_\xi)) \\
  &=& \delta_{i,h_a \cdot j} \1z_i \otimes ((e_i U_{h_a})(d_\xi U_{h_\xi})(h)) \\
  &=& \rukl{\Theta^{-1}((\1z_i \sigma_{h_a} \otimes e_i U_{h_a}) (\1z_j \otimes (d_\xi U_{h_\xi}))}(h).
\eglnoz
Thus $\Theta \circ \psi(a) \circ \Theta^{-1} = (\1z_i \otimes e_i) \circ (\sigma_{h_a} \otimes U_{h_a})$.

So all in all, we have computed that $j_r^H(\res_H^G(\xf)) \in KK(c_0(I) \rtimes_{\sigma,r} H, D \rtimes_{\tau,r} H)$ is represented by the Kasparov module
\bgloz
  (\ell^2(I,(D \rtimes_{\tau,r} H)), \phi \rtimes_r H, 0)
\egloz
where $\phi \rtimes_r H: c_0(I) \rtimes_{\sigma,r} H \to \cL(\ell^2(I,(D \rtimes_{\tau,r} H)))$ is given by
\bgloz
  c_0(I) \rtimes_{\sigma,r} H \ni \1z_i U_h \ma (\1z_i \otimes e_i) \circ (\sigma_h \otimes U_h) \in \cL(\ell^2(I,(D \rtimes_{\tau,r} H))).
\egloz
\eproof

\subsection{Direct sum decomposition}
\label{dirsumdec}

Let $H$ be a subgroup of $G$. It is possible to decompose $c_0(I) \rtimes_{\sigma,r} H$ into direct summands corresponding to the $H$-orbits on $I$, i.e.
\bgloz
  c_0(I) \rtimes_{\sigma,r} H = \bigoplus_{\eckl{i} \in H \backslash I} \rukl{c_0(H \cdot i) \rtimes_{\sigma,r} H}.
\egloz
Let us denote the summand $c_0(H \cdot i) \rtimes_{\sigma,r} H$ corresponding to $\eckl{i} \in H \backslash I$ by $C_{\eckl{i}}$ and let $\iota_{\eckl{i}}$ be the embedding $C_{\eckl{i}} \to c_0(I) \rtimes_{\sigma,r} H$.

As explained in \cite{Bla}, Theorem~19.7.1,
\bgloz
  \prod_{\eckl{i} \in H \backslash I} (KK(\iota_{\eckl{i}}) \otimes \sqcup): 
  KK(c_0(I) \rtimes_{\sigma,r} H, D \rtimes_{\tau,r} H) \to \prod_{\eckl{i} \in H \backslash I} KK(C_{\eckl{i}}, D \rtimes_{\tau,r} H)
\egloz
is an isomorphism. Here \an{$\otimes$} stands for the Kasparov product.

It is immediate that under this isomorphism, the element $j_r^H(\res_H^G(\xf))$ corresponds to $(\xf_{\eckl{i}})_{\eckl{i} \in H \backslash I}$ where $\xf_{\eckl{i}} \in KK(C_{\eckl{i}}, D \rtimes_{\tau,r} H)$ is represented by the Kasparov module
\bgl
\label{x_i}
  (\ell^2(H \cdot i,(D \rtimes_{\tau,r} H)), (\phi \rtimes_r H)_{\eckl{i}}, 0)
\egl
with $(\phi \rtimes_r H)_{\eckl{i}}$ given by $C_{\eckl{i}} \ni \1z_j U_h \ma (\1z_j \otimes e_j) \circ (\sigma_h \otimes U_h) \in \cL(\ell^2(H \cdot i,(D \rtimes_{\tau,r} H)))$.

In other words, we have
\bgl
\label{x_i=iota_ix}
  \xf_{\eckl{i}} = KK(\iota_{\eckl{i}}) \otimes j_r^H(\res_H^G(\xf)).
\egl
We describe $\xf_{\eckl{i}}$ alternatively as follows: Let $\varphi_{\eckl{i}}$ be the homomorphism
\bgln
\label{phi_[i]}
  \varphi_{\eckl{i}}: C_{\eckl{i}} &\to& \cK(\ell^2(H \cdot i)) \otimes (D \rtimes_{\tau,r} H) \\
  \1z_j U_h &\ma& e_{j,h^{-1} \cdot j} \otimes e_j U_h \nonumber
\egln
where $e_{j,h^{-1} \cdot j}$ is the rank $1$ operator $\spkl{\sqcup, \varepsilon_{h^{-1} \cdot j}} \varepsilon_j \in \cL(\ell^2(H \cdot i))$ ($\menge{\varepsilon_j}{j \in H \cdot i}$ is the canonical orthonormal basis of $\ell^2(H \cdot i)$).

Existence of $\varphi_{\eckl{i}}$ can be seen as follows: Using a faithful representation of $D$ on a Hilbert space $\cH$, we can view $D$ as a sub-C*-algebra of $\cL(\cH)$. Hence, according to the definition of the reduced crossed product, the C*-algebra $\cK(\ell^2(H \cdot i)) \otimes_{\min} (D \rtimes_{\tau,r} G)$ acts on the Hilbert space $\ell^2(H \cdot i) \otimes \cH \otimes \ell^2(H)$. At the same time, using the definition of the reduced crossed product $C_{\eckl{i}} = c_0(H \cdot i) \rtimes_{\sigma,r} H$, we obtain a faithful representation $\pi$ of $C_{\eckl{i}}$ sending $\1z_j U_h \in C_{\eckl{i}}$ to the operator $\pi(\1z_j) (1 \otimes 1 \otimes \lambda_h)$ on $\ell^2(H \cdot i) \otimes \cH \otimes \ell^2(H)$ where $\pi(\1z_j)$ is given by $\pi(\1z_j)(\ve_k \otimes \xi \otimes \ve_x) = (e_{x^{-1} \cdot j, x^{-1} \cdot j} \otimes e_{x^{-1} \cdot j} \otimes 1)(\ve_k \otimes \xi \otimes \ve_x)$ for $j, k \in H \cdot i$, $\xi \in \cH$ and $x \in H$. Here $\gekl{\ve_{h \cdot i}}_{h \in H}$ is the canonical orthonormal basis of $\ell^2(H \cdot i)$, $\gekl{\ve_x}_{x \in H}$ is the canonical orthonormal basis of $\ell^2(H)$ and $e_{x^{-1} \cdot j, x^{-1} \cdot j}$ is the rank $1$ projection corresponding to the basis vector $\ve_{x^{-1} \cdot j} \in \ell^2(H \cdot i)$. Now, applying Fell's absorption principle or rather adapting its proof, we consider the unitary $W$ on $\ell^2(H \cdot i) \otimes \cH \otimes \ell^2(H)$ given by $W(\ve_k \otimes \xi \otimes \ve_x) = \ve_{x \cdot k} \otimes \xi \otimes \ve_x$. Then a direct computation shows
\bgloz
  \Ad(W) \circ \rukl{\pi(\1z_j) (1 \otimes 1 \otimes \lambda_h)} = e_{j, h^{-1} \cdot j} \otimes e_j U_h.
\egloz
Therefore, $\Ad(W) \circ \pi$ is the desired homomorphism $\varphi_{\eckl{i}}$.

The homomorphism $e_{i,i} \otimes \id_{D \rtimes_{\tau,r} H}: D \rtimes_{\tau,r} H \to \cK(\ell^2(H \cdot i)) \otimes (D \rtimes_{\tau,r} H)$; $b \ma e_{i,i} \otimes b$ gives a KK-equivalence between $D \rtimes_{\tau,r} H$ and $\cK(\ell^2(H \cdot i)) \otimes (D \rtimes_{\tau,r} H)$.

\bsublemma
\label{desc-x_i}
$\xf_{\eckl{i}} = KK(\varphi_{\eckl{i}}) \otimes KK(e_{i,i} \otimes \id_{D \rtimes_{\tau,r} H})^{-1}$ where $\otimes$ is the Kasparov product.
\esublemma
\bproof
Viewing $D \rtimes_{\tau,r} H$ as a full corner in $\cK(\ell^2(H \cdot i)) \otimes (D \rtimes_{\tau,r} H)$ via $e_{i,i} \otimes \id_{D \rtimes_{\tau,r} H}$, it is clear that $KK(e_{i,i} \otimes \id_{D \rtimes_{\tau,r} H})^{-1}$ is represented by the Kasparov module given by the $\rukl{\cK(\ell^2(H \cdot i)) \otimes (D \rtimes_{\tau,r} H)}$--$D \rtimes_{\tau,r} H$-imprimitivity bimodule $\ell^2(H \cdot i,(D \rtimes_{\tau,r} H))$. This Kasparov module is explicitly given by the Hilbert $D \rtimes_{\tau,r} H$-module $\ell^2(H \cdot i,(D \rtimes_{\tau,r} H))$ and the left action
\bgloz
  \cK(\ell^2(H \cdot i)) \otimes (D \rtimes_{\tau,r} H) \to \cK(\ell^2(H \cdot i,(D \rtimes_{\tau,r} H))); \: e_{j,h^{-1} \cdot j} \otimes b \ma \1z_j \sigma_h \otimes b.
\egloz
Using the descriptions of $\xf_{\eckl{i}}$ and $\varphi_{\eckl{i}}$ from \eqref{x_i} and \eqref{phi_[i]}, it is clear that
\bgloz
  \xf_{\eckl{i}} = KK(\varphi_{\eckl{i}}) \otimes KK(e_{i,i} \otimes \id_{D \rtimes_{\tau,r} H})^{-1}.
\egloz
\eproof

\bsubcor
\label{phi_i^B}
Let $B$ be a sub-C*-algebra of $D \rtimes_{\tau,r} H$ such that for all $j \in H \cdot i$ and $h \in H$, $e_j U_h$ lies in $B$. Let $\iota$ be the inclusion $B \into D \rtimes_{\tau,r} H$, let $\varphi_{\eckl{i}} \vert^B$ be the homomorphism $C_{\eckl{i}} \to \cK(\ell^2(H \cdot i)) \otimes B$; $a \ma \varphi_{\eckl{i}}(a)$ (we just restrict the image of $\varphi_{\eckl{i}}$) and denote by $e_{i,i} \otimes \id_B$ the homomorphism $B \to \cK(\ell^2(H \cdot i)) \otimes B$; $b \ma e_{i,i} \otimes b$. Then
\bgloz
  KK(\varphi_{\eckl{i}} \vert^B) \otimes KK(e_{i,i} \otimes \id_B)^{-1} \otimes KK(\iota) = \xf_{\eckl{i}}.
\egloz
\esubcor
\bproof
We have
\bglnoz
  && KK(\varphi_{\eckl{i}} \vert^B) \otimes KK(e_{i,i} \otimes \id_B)^{-1} \otimes KK(\iota) \otimes KK(e_{i,i} \otimes \id_{D \rtimes_{\tau,r} H}) \\
  &=& KK(\varphi_{\eckl{i}} \vert^B) \otimes KK(e_{i,i} \otimes \id_B)^{-1} \otimes KK(e_{i,i} \otimes \id_B) \otimes KK(\id_{\cK(\ell^2(H \cdot i))} \otimes \iota) \\
  &=& KK(\varphi_{\eckl{i}} \vert^B) \otimes KK(\id_{\cK(\ell^2(H \cdot i))} \otimes \iota) \\
  &=& KK(\varphi_{\eckl{i}}).
\eglnoz
Now multiply with $KK(e_{i,i} \otimes \id_{D \rtimes_{\tau,r} H})^{-1}$ from the right and use Lemma~\ref{desc-x_i}.
\eproof

\subsection{KK-equivalences for all finite subgroups}
\label{proofmainthm-end}

Now we consider finite subgroups. Since in this case, we do not have to distinguish between full and reduced crossed products, we can omit the index $r$ everywhere. Our goal is to prove 
\bsubtheo
For every finite subgroup $H$ of $G$, the element $j^H(\res_H^G(\xf))$ in $KK(c_0(I) \rtimes_{\sigma} H, D \rtimes_{\tau} H)$ is a KK-equivalence.
\esubtheo
As both $c_0(I) \rtimes_\sigma H$ and $D \rtimes_\tau H$ satisfy the UCT being crossed products of commutative C*-algebras by amenable groups, it suffices to prove that $j^H(\res_H^G(\xf))$ induces an isomorphism on K-theory. To show this, the strategy is to reduce everything to finite dimensional sub-C*-algebras. Therefore, we write both $c_0(I) \rtimes_\sigma H$ and $D \rtimes_\tau H$ as inductive limits of finite dimensional C*-algebras and consider the corresponding inductive limit descriptions of their K-groups.

In the sequel, we write $K_*$ for the direct sum of $K_0$ of $K_1$ viewed as a $\Zz / 2 \Zz$-graded abelian group.

We start with $c_0(I) \rtimes_\sigma H$. We have already seen in \S~\ref{dirsumdec} the decomposition $c_0(I) \rtimes_\sigma H = \bigoplus_{\eckl{i} \in H \backslash I} C_{\eckl{i}}$. Thus, it is clear that we have $c_0(I) \rtimes_\sigma H \cong \ilim_F \bigoplus_{\eckl{i} \in \eckl{F}} C_{\eckl{i}}$, where the limit is taken over the finite subsets $F$ of $I$ and we denote the image of $F$ under the projection $I \to H \backslash I$ by $\eckl{F}$. Therefore we obtain $\ilim_F \bigoplus_{\eckl{i} \in \eckl{F}} K_*(C_{\eckl{i}}) \cong K_*(c_0(I) \rtimes_\sigma H)$, and this identification is induced by the homomorphisms
\bgloz
  \sum_{\eckl{i} \in \eckl{F}} \iota_{\eckl{i}}: \bigoplus_{\eckl{i} \in \eckl{F}} K_*(C_{\eckl{i}}) \to K_*(c_0(I) \rtimes_\sigma H).
\egloz

Now we consider $D \rtimes_\tau H$. For a finite subset $F$ of $I$, let $(D \rtimes_\tau H)_F$ be the sub-C*-algebra of $D \rtimes_\tau H$ which is generated by $\menge{e_i U_h}{i \in H \cdot F, h \in H}$. As before, we certainly have $D \rtimes_\tau H \cong \ilim_F (D \rtimes_\tau H)_F$ and thus $\ilim_F K_*((D \rtimes_\tau H)_F) \cong K_*(D \rtimes_\tau H)$. This identification is realized by the homomorphisms induced by the canonical inclusions $(D \rtimes_\tau H)_F \into D \rtimes_\tau H$ on K-theory.

We now compare these direct limit decompositions. Given a finite subset $F$ of $I$, we set
\bgl
\label{xf_i^F}
  \xf_{\eckl{i}}^F \defeq KK(\varphi_{\eckl{i}} \vert^{(D \rtimes_\tau H)_F}) \otimes KK(e_{i,i} \otimes \id_{(D \rtimes_\tau H)_F})^{-1}
\egl
using the notation from Corollary~\ref{phi_i^B}. Let $K_*(\xf_{\eckl{i}}^F)$ be the homomorphism induced on K-theory by $\xf_{\eckl{i}}^F$. By \eqref{x_i=iota_ix} and Corollary~\ref{phi_i^B}, the diagram
\bgl
\label{cd}
  \begin{CD}
  K_*(C_{\eckl{i}}) @> K_*(\xf_{\eckl{i}}^F) >> K_*((D \rtimes_\tau H)_F) \\
  @VV K_*(\iota_{\eckl{i}}) V @VVV \\
  K_*(c_0(I) \rtimes_\sigma H) @> K_*(j^H(\res_H^G(\xf))) >> K_*(D \rtimes_\tau H)
  \end{CD}
\egl
commutes, where the right vertical arrow is induced by the canonical inclusion $(D \rtimes_\tau H)_F \into D \rtimes_\tau H$. Therefore, for every finite subset $F$ of $I$, we have a homomorphism $\sum_{\eckl{i} \in \eckl{F}} K_*(\xf_{\eckl{i}}^F): K_*(\bigoplus_{\eckl{i} \in \eckl{F}} C_{\eckl{i}}) \to K_*((D \rtimes_\tau H)_F)$, and these homomorphisms induce a homomorphism
\bgloz
  \ilim_F \sum_{\eckl{i} \in \eckl{F}} K_*(\xf_{\eckl{i}}^F): \ilim_F \bigoplus_{\eckl{i} \in \eckl{F}} K_*(C_{\eckl{i}}) \to \ilim_F K_*((D \rtimes_\tau H)_F)
\egloz
by a similar computation as in Corollary~\ref{phi_i^B}. By commutativity of \eqref{cd}, the diagram
\bgl
\label{cd2}
  \begin{CD}
  \ilim_F \bigoplus_{\eckl{i} \in \eckl{F}} K_*(C_{\eckl{i}}) @> \ilim_F \sum_{\eckl{i} \in \eckl{F}} K_*(\xf_{\eckl{i}}^F) >> \ilim_F K_*((D \rtimes_\tau H)_F) \\
  @VV \cong V @VV \cong V \\
  K_*(c_0(I) \rtimes_\sigma H) @> K_*(j^H(\res_H^G(\xf))) >> K_*(D \rtimes_\tau H)
  \end{CD}
\egl
commutes as well.

In these inductive limits, it clearly suffices to only take those finite subsets $F$ which satisfy the condition that $\menge{e_i}{i \in H \cdot F} \cup \gekl{0}$ is multiplicatively closed. Now the point is that we will prove in the next proposition that for these finite subsets $F$, $\sum_{\eckl{i} \in \eckl{F}} K_*(\xf_{\eckl{i}}^F): \bigoplus_{\eckl{i} \in \eckl{F}} K_*(C_{\eckl{i}}) \to K_*((D \rtimes_\tau H)_F)$ is an isomorphism. This will then imply that the homomorphism $\ilim_F \sum_{\eckl{i} \in \eckl{F}} K_*(\xf_{\eckl{i}}^F): \ilim_F \bigoplus_{\eckl{i} \in \eckl{F}} K_*(C_{\eckl{i}}) \to \ilim_F K_*((D \rtimes_\tau H)_F)$ is an isomorphism, where we take the inductive limit over those $F$ satisfying the condition that $\menge{e_i}{i \in H \cdot F} \cup \gekl{0}$ is multiplicatively closed. Because diagram~\eqref{cd2} commutes, this will then imply our main observation that $K_*(j^H(\res_H^G(\xf)))$ is an isomorphism.

\bsubprop
\label{KKequi-finite}
Let $F$ be a finite subset of $I$ such that $\menge{e_j}{j \in H \cdot F} \cup \gekl{0}$ is multiplicatively closed. Then the KK-elements $\xf_{\eckl{i}}^F$, $\eckl{i} \in \eckl{F}$, induce a K-theoretic isomorphism
\bgloz
  \sum_{\eckl{i} \in \eckl{F}} K_*(\xf_{\eckl{i}}^F): \bigoplus_{\eckl{i} \in \eckl{F}} K_*(C_{\eckl{i}}) \to K_*((D \rtimes_\tau H)_F).
\egloz
\esubprop
\bproof
We decompose $(D \rtimes_\tau H)_F$ into direct summands as follows: For every $j \in F$, set $e(j) \defeq e_j - \bigvee_{k \in H \cdot F, e_k \lneq e_j} e_k$ and for every $i \in F$, define $e(\eckl{i}) = \sum_{j \in H \cdot i} e(j)$. By construction, the following facts hold:
\begin{itemize}
\item For every $j$ in $F$, $e(j) \neq 0$ as $\gekl{e_i}_{i \in I}$ is independent (see (II)) and because of our assumption that $e_j \neq 0$ for all $j \in I$.
\item For $i,j \in F$ with $\eckl{i} \neq \eckl{j}$, $e(\eckl{i}) \perp e(\eckl{j})$.
\item $\bigvee_{j \in H \cdot F} e_j = \sum_{\eckl{i} \in \eckl{F}} e(\eckl{i})$.
\item For every $i$ in $F$, $e(\eckl{i})$ is $H$-invariant with respect to the action $\tau$.
\end{itemize}

The last fact implies that these projections $e(\eckl{i})$ are central in $(D \rtimes_\tau H)_F$. Thus, using this, the second and third fact and also our condition that $\menge{e_i}{i \in H \cdot F} \cup \gekl{0}$ is multiplicatively closed, we deduce
\bgloz
  (D \rtimes_\tau H)_F = \bigoplus_{\eckl{i} \in \eckl{F}} e(\eckl{i}) ((D \rtimes_\tau H)_F) e(\eckl{i}).
\egloz
Using the first two facts, it is straightforward to check that $e(\eckl{i}) ((D \rtimes_\tau H)_F) e(\eckl{i})$ is generated as a C*-algebra by the elements $e(j) U_h$ for $j \in H \cdot i$, $h \in H$, and that we can identify $e(\eckl{i}) ((D \rtimes_\tau H)_F) e(\eckl{i})$ with $c_0(H \cdot i) \rtimes_\sigma H = C_{\eckl{i}}$ via
\bgl
\label{pi_F^[i]}
  e(j) U_h \ma \1z_j U_h.
\egl
Thus we obtain an isomorphism
\bgloz
  (D \rtimes_\tau H)_F \cong \bigoplus_{\eckl{i} \in \eckl{F}} C_{\eckl{i}}.
\egloz
Let $\pi_F^{\eckl{i}}: (D \rtimes_\tau H)_F \cong \bigoplus_{\eckl{i} \in \eckl{F}} C_{\eckl{i}} \to C_{\eckl{i}}$ be the composition of this isomorphism with the canonical projection $\bigoplus_{\eckl{i} \in \eckl{F}} C_{\eckl{i}} \to C_{\eckl{i}}$. It follows that in K-theory, $\bigoplus_{\eckl{i} \in \eckl{F}} K_*(\pi_F^{\eckl{i}}): K_*((D \rtimes_\tau H)_F) \to \bigoplus_{\eckl{i} \in \eckl{F}} K_*(C_{\eckl{i}})$ is an isomorphism. This means that to show that $\sum_{\eckl{i} \in \eckl{F}} K_*(\xf_{\eckl{i}}^F)$ is an isomorphism, we can equally well prove that the composition
\bgloz
  \bigoplus_{\eckl{i} \in \eckl{F}} K_*(C_{\eckl{i}}) \overset{\sum_{\eckl{i} \in \eckl{F}} K_*(\xf_{\eckl{i}}^F)}{\lori} K_*((D \rtimes_\tau H)_F)
  \overset{\bigoplus_{\eckl{i} \in \eckl{F}} K_*(\pi_F^{\eckl{i}})}{\lori} \bigoplus_{\eckl{i} \in \eckl{F}} K_*(C_{\eckl{i}})
\egloz
is an isomorphism.

This composition can be described by a $\eckl{F} \times \eckl{F}$-matrix whose $(\eckl{i},\eckl{j})$-th entry is given by $K_*(\pi_F^{\eckl{i}}) \circ K_*(\xf_{\eckl{j}}^F)$ (here $\circ$ is composition of homomorphisms). Going through our constructions, it is clear that 
\bgl
\label{uppertri}
K_*(\pi_F^{\eckl{i}}) \circ K_*(\xf_{\eckl{j}}^F) \neq 0 \text{ only if } \bigvee_{l \in H \cdot j} e_l \geq e(\eckl{i}) \LRarr \bigvee_{l \in H \cdot j} e_l \geq \bigvee_{k \in H \cdot i} e_k.
\egl
It is immediate that $\eckl{j} \geq \eckl{i} :\LRarr \bigvee_{l \in H \cdot j} e_l \geq \bigvee_{k \in H \cdot i} e_k$ defines a partial order relation on $\eckl{F}$. If we arrange the elements of $\eckl{F}$ in increasing order with respect to this partial order (increasing means that the elements $\eckl{j}$ which come after an element $\eckl{i}$ do not satisfy $\eckl{j} \leq \eckl{i}$), then \eqref{uppertri} tells us that $(K_*(\pi_F^{\eckl{i}}) \circ K_*(\xf_{\eckl{j}}^F))_{\eckl{i},\eckl{j}}$ becomes an upper triangular matrix. Hence the $\eckl{F} \times \eckl{F}$-matrix describing $\rukl{\bigoplus_{\eckl{i} \in \eckl{F}} K_*(\pi_F^{\eckl{i}})} \circ \rukl{\sum_{\eckl{i} \in \eckl{F}} K_*(\xf_{\eckl{i}}^F)}$ is the sum of a nilpotent matrix and a diagonal matrix whose $(\eckl{i},\eckl{i})$-th entry is given by $K_*(\pi_F^{\eckl{i}}) \circ K_*(\xf_{\eckl{i}}^F)$. To prove that the matrix $\rukl{K_*(\pi_F^{\eckl{i}}) \circ K_*(\xf_{\eckl{j}}^F)}_{\eckl{i},\eckl{j}}$ is invertible, it remains to prove that the diagonal entries of this matrix are invertible, i.e. that $K_*(\pi_F^{\eckl{i}}) \circ K_*(\xf_{\eckl{i}}^F): K_*(C_{\eckl{i}}) \to K_*(C_{\eckl{i}})$ is an isomorphism for all $\eckl{i} \in \eckl{F}$. 

Recall that
\bgloz
  \xf_{\eckl{i}}^F = KK(\varphi_{\eckl{i}} \vert^{(D \rtimes_\tau H)_F}) \otimes KK(e_{i,i} \otimes \id_{(D \rtimes_\tau H)_F})^{-1}
\egloz
so that
\bgloz
  K_*(\xf_{\eckl{i}}^F) = K_*(e_{i,i} \otimes \id_{(D \rtimes_\tau H)_F})^{-1} \circ K_*(\varphi_{\eckl{i}} \vert^{(D \rtimes_\tau H)_F}).
\egloz
As
\bgloz
  (e_{i,i} \otimes \id_{C_{\eckl{i}}}) \circ \pi_F^{\eckl{i}} = (\id_{\cL(\ell^2(H \cdot i))} \otimes \pi_F^{\eckl{i}}) \circ (e_{i,i} \otimes \id_{(D \rtimes_\tau H)_F}),
\egloz
we obtain
\bglnoz
  && K_*(\pi_F^{\eckl{i}}) \circ K_*(\xf_{\eckl{i}}^F) \\
  &=& K_*(\pi_F^{\eckl{i}}) \circ K_*(e_{i,i} \otimes \id_{(D \rtimes_\tau H)_F})^{-1} \circ K_*(\varphi_{\eckl{i}} \vert^{(D \rtimes_\tau H)_F}) \\
  &=& K_*((e_{i,i} \otimes \id_{C_{\eckl{i}}}))^{-1} \circ K_*(\id_{\cL(\ell^2(H \cdot i))} \otimes \pi_F^{\eckl{i}}) 
  \circ K_*(\varphi_{\eckl{i}} \vert^{(D \rtimes_\tau H)_F}) \\
  &=& K_*((e_{i,i} \otimes \id_{C_{\eckl{i}}}))^{-1} 
  \circ K_*((\id_{\cL(\ell^2(H \cdot i))} \otimes \pi_F^{\eckl{i}}) \circ \varphi_{\eckl{i}} \vert^{(D \rtimes_\tau H)_F}).
\eglnoz
Note that $\cL(\ell^2(H \cdot i)) = \cK(\ell^2(H \cdot i))$ as $H$ is finite.

To prove that $K_*(\pi_F^{\eckl{i}}) \circ K_*(\xf_{\eckl{i}}^F)$ is an isomorphism, it therefore suffices to check that $K_*((\id_{\cL(\ell^2(H \cdot i))} \otimes \pi_F^{\eckl{i}}) \circ \varphi_{\eckl{i}} \vert^{(D \rtimes_\tau H)_F})$ is an isomorphism. By \eqref{phi_[i]} and \eqref{pi_F^[i]}, $(\id_{\cL(\ell^2(H \cdot i))} \otimes \pi_F^{\eckl{i}}) \circ \varphi_{\eckl{i}} \vert^{(D \rtimes_\tau H)_F}$ is given by the homomorphism
\bgloz
  C_{\eckl{i}} \to \cL(\ell^2(H \cdot i)) \otimes C_{\eckl{i}}; \: \1z_j U_h \ma e_{j,h^{-1} \cdot j} \otimes \1z_j U_h.
\egloz
Let $s: H \cdot i \to H$ be a map satisfying $s(h \cdot i) \cdot i = h \cdot i$. Define
\bgloz
  W \defeq \sum_{j \in H \cdot i} \sigma_{s(j)} \otimes \1z_j \in \cL(\ell^2(H \cdot i)) \otimes C_{\eckl{i}}. 
\egloz
We finally claim that $W$ is a unitary such that
\bgloz
  \Ad(W^*) \circ (\id_{\cL(\ell^2(H \cdot i))} \otimes \pi_F^{\eckl{i}}) \circ \varphi_{\eckl{i}} \vert^{(D \rtimes_\tau H)_F} = e_{i,i} \otimes \id_{C_{\eckl{i}}}.
\egloz
This follows from the following computations:
\bglnoz
  W^* W &=& \sum_{j \in H \cdot i} \sigma_{s(j)}^* \sigma_{s(j)} \otimes \1z_j = 1 \otimes 1, \\
  W W^* &=& \sum_{j \in H \cdot i} \sigma_{s(j)} \sigma_{s(j)}^* \otimes \1z_j = 1 \otimes 1
\eglnoz
and
\bglnoz
  && W^* \rukl{e_{j,h^{-1} \cdot j} \otimes \1z_j U_h} W 
  = \sum_{k,l \in H \cdot i} \sigma_{s(k)}^* e_{j,h^{-1} \cdot j} \sigma_{s(l)} \otimes \underbrace{\1z_k \1z_j U_h \1z_l}_{= \delta_{k,j} \delta_{j,hl} \1z_j U_h} \\
  &=& \sigma_{s(j)}^* e_{j,h^{-1} \cdot j} \sigma_{s(h^{-1} \cdot j)} \otimes \1z_j U_h = e_{i,i} \otimes \1z_j U_h = (e_{i,i} \otimes \id_{C_{\eckl{i}}})(\1z_j U_h).
\eglnoz
Thus $(\id_{\cL(\ell^2(H \cdot i))} \otimes \pi_F^{\eckl{i}}) \circ \varphi_{\eckl{i}} \vert^{(D \rtimes_\tau H)_F}$ is unitarily equivalent to $e_{i,i} \otimes \id_{C_{\eckl{i}}}$. As $e_{i,i} \otimes \id_{C_{\eckl{i}}}$ induces an isomorphism on K-theory, we are done.
\eproof

\bproof[Proof of Theorem~\ref{mainthm}]
Theorem~\ref{mainthm} now follows from Proposition~\ref{KKequi-finite} and commutativity of diagram~\eqref{cd2}.
\eproof

\subsection{Baum-Connes}
\label{B-C}

Under certain K-theoretic assumptions on our group $G$, we may now apply the Baum-Connes machinery to our situation.
\bsubcor
\label{cor-mainthm1}
Let us assume that the conditions of Theorem~\ref{mainthm} are satisfied, i.e. that conditions (I), (II) and (III) from \S~\ref{generalK} hold. Moreover, assume that the group $G$ satisfies the Baum-Connes conjecture with coefficients in $c_0(I)$ and $D$ with respect to the $G$-actions $\sigma$ and $\tau$. Then the descent $j_r^G(\xf) \in KK(c_0(I) \rtimes_{\sigma,r} G, D \rtimes_{\tau,r} G)$ induces an isomorphism on K-theory.
\esubcor
\bproof
We have proven in Theorem~\ref{mainthm} that for all finite subgroups $H$ of $G$, the descent $j^H(\res_H^G(\xf))$ is a KK-equivalence. Now our corollary follows from \cite{E-N-O}, Proposition~2.1~(i).
\eproof
Under additional assumptions, we even obtain
\bsubcor
\label{cor-mainthm2}
If we assume, in addition to the requirements of the previous corollary, that both reduced crossed products $c_0(I) \rtimes_{\sigma,r} G$ and $D \rtimes_{\tau,r} G$ satisfy the UCT, then $j_r^G(\xf)$ is a KK-equivalence.
\esubcor
\bproof
This follows immediately from the previous corollary.
\eproof
To obtain that $j_r^G(\xf)$ is a KK-equivalence, we can also proceed as follows:
\bsubcor
\label{cor-mainthm3}
If we assume, in addition to the requirements of Corollary~\ref{cor-mainthm1}, that $G$ satisfies the strong Baum-Connes conjecture with coefficients in $c_0(I)$ and $D$ with respect to the $G$-actions $\sigma$ and $\tau$, then $j_r^G(\xf)$ is a KK-equivalence.
\esubcor
\bproof
This follows from Theorem~\ref{mainthm} and \cite{E-N-O}, Proposition~2.1~(iii).
\eproof
The conditions of this corollary are for instance satisfied if $G$ is amenable.

\subsection{Imprimitivity Theorems}
\label{impthm}

Consider the direct sum decomposition
\bgloz
  c_0(I) \rtimes_{\sigma,r} G = \bigoplus_{\eckl{i} \in G \backslash I} \rukl{c_0(G \cdot i) \rtimes_{\sigma,r} G}.
\egloz
As before, we denote the summand $c_0(G \cdot i) \rtimes_{\sigma,r} G$ corresponding to $\eckl{i} \in G \backslash I$ by $C_{\eckl{i}}$ and let $\iota_{\eckl{i}}$ be the embedding $C_{\eckl{i}} \to c_0(I) \rtimes_{\sigma,r} G$. 

Under the isomorphism
\bgloz
  \prod_{\eckl{i} \in G \backslash I} (KK(\iota_{\eckl{i}}) \otimes \sqcup): 
  KK(c_0(I) \rtimes_{\sigma,r} G, D \rtimes_{\tau,r} G) \to \prod_{\eckl{i} \in G \backslash I} KK(C_{\eckl{i}}, D \rtimes_{\tau,r} G)
\egloz
from \cite{Bla}, Theorem~19.7.1, the element $j_r^G(\xf) \in KK(c_0(I) \rtimes_{\sigma,r} G, D \rtimes_{\tau,r} G)$ corresponds to $(\xf_{\eckl{i}})_{\eckl{i} \in G \backslash I}$ with $\xf_{\eckl{i}} \defeq KK(\iota_{\eckl{i}}) \otimes j_r^G(\xf)$. By Lemma~\ref{desc-x_i}, we have
\bgl
\label{x_i-for-G}
  \xf_{\eckl{i}} = KK(\varphi_{\eckl{i}}) \otimes KK(e_{i,i} \otimes \id_{D \rtimes_{\tau,r} G})^{-1}
\egl
where the homomorphisms $\varphi_{\eckl{i}}$ and $e_{i,i} \otimes \id_{D \rtimes_{\tau,r} G}$ are given by
\bgloz
  \varphi_{\eckl{i}}: C_{\eckl{i}} \to \cK(\ell^2(G \cdot i)) \otimes_{\min} (D \rtimes_{\tau,r} G), \: \1z_j U_g \ma e_{j,g^{-1} \cdot j} \otimes e_j U_g
\egloz
and
\bgloz
e_{i,i} \otimes \id_{D \rtimes_{\tau,r} G}: D \rtimes_{\tau,r} G \to \cK(\ell^2(G \cdot i)) \otimes_{\min} (D \rtimes_{\tau,r} G), \: T \ma e_{i,i} \otimes T.
\egloz

To further examine $\xf_{\eckl{i}}$, let us now describe $C_{\eckl{i}} = c_0(G \cdot i) \rtimes_{\sigma,r} G$ up to Morita equivalence with the help of concrete homomorphisms. For $i \in I$, let $G_i$ be the stabilizer of $i$, i.e. $G_i \defeq \menge{g \in G}{g \cdot i = i}$. Then we have a bijection $G / G_i \cong G \cdot i$, $g G_i \ma g \cdot i$ which is $G$-equivariant. Thus we can identify $C_{\eckl{i}} = c_0(G \cdot i) \rtimes_{\sigma,r} G$ with $c_0(G / G_i) \rtimes_r G$ where we take the translation action of $G$ on $c_0(G / G_i)$ for the second reduced crossed product. Moreover, the homomorphism
\bgloz
  C^*_r(G_i) \to c_0(G / G_i) \rtimes_r G; \: \lambda_g \ma \1z_{eG_i} U_g
\egloz
exists by Lemma~\ref{functor} and induces a KK-equivalence in $KK(C^*_r(G_i), c_0(G / G_i) \rtimes_r G)$. The last assertion follows from the observation that the projection $\1z_{e G_i} \in c_0(G / G_i) \rtimes_r G$ is a full projection and that the above homomorphism yields an isomorphism $C^*_r(G_i) \cong \1z_{e G_i} \rukl{c_0(G / G_i) \rtimes_r G} \1z_{e G_i}; \lambda_g \ma \1z_{eG_i} U_g$ (injectivity follows from Lemma~\ref{functor} and surjectivity can be seen immediately).

Composing the homomorphism $C^*_r(G_i) \to c_0(G / G_i) \rtimes_r G$ and the canonical identification $c_0(G / G_i) \rtimes_r G \cong c_0(G \cdot i) \rtimes_{\sigma,r} G = C_{\eckl{i}}$ from above, we obtain the homomorphism
\bgl
\label{phii}
  \varphi_i: C^*_r(G_i) \to C_{\eckl{i}}, \: \lambda_g \ma \1z_i U_g.
\egl
By our observations, $KK(\varphi_i)$ is a KK-equivalence in $KK(C^*_r(G_i),C_{\eckl{i}})$. For two different choices of the representative $i$ of the class $\eckl{i} \in G \backslash I$, the stabilizers will be different in general, but they will always be conjugate. So the choices of the particular representatives do not really matter.

Finally, let us compute the Kasparov product $KK(\varphi_i) \otimes \xf_{\eckl{i}}$. As a preparation, note that $(\varphi_{\eckl{i}} \circ \varphi_i) (\lambda_g) = e_{i,i} \otimes e_i U_g \in \cK(\ell^2(G \cdot i)) \otimes_{\min} (D \rtimes_{\tau,r} G)$ for $g \in G_i$. Thus composing $\varphi_{\eckl{i}} \circ \varphi_i$ with the canonical identification $e_{i,i} \otimes D \rtimes_{\tau,r} G \cong D \rtimes_{\tau,r} G$, we obtain a homomorphism $\Phi_i: C^*_r(G_i) \to D \rtimes_{\tau,r} G, \lambda_g \ma e_i U_g$. By construction,
\bgl
\label{phiphi}
  \varphi_{\eckl{i}} \circ \varphi_i = (e_{i,i} \otimes \id_{D \rtimes_{\tau,r} G}) \circ \Phi_i.
\egl
Thus
\bgln
\label{phix_i}
  && KK(\varphi_i) \otimes \xf_{\eckl{i}} = KK(\varphi_i) \otimes KK(\iota_{\eckl{i}}) \otimes j_r^G(\xf) \\
  &\overset{\eqref{x_i-for-G}}{=}& KK(\varphi_i) \otimes KK(\varphi_{\eckl{i}}) \otimes KK(e_{i,i} \otimes \id_{D \rtimes_{\tau,r} G})^{-1} 
  \nonumber \\
  &\overset{\eqref{phiphi}}{=}& KK(\Phi_i) \otimes KK(e_{i,i} \otimes \id_{D \rtimes_{\tau,r} G}) \otimes KK(e_{i,i} \otimes \id_{D \rtimes_{\tau,r} G})^{-1} \nonumber \\
  &=& KK(\Phi_i). \nonumber
\egln

Let us summarize our observations.
\bsubprop
\label{imp}
Let $\cR$ be a complete system of representatives for $G \backslash I$. The homomorphism
\bgloz
  \bigoplus_{i \in \cR} \varphi_i: \bigoplus_{i \in \cR} C^*_r(G_i) \to \bigoplus_{i \in \cR} C_{\eckl{i}} = c_0(I) \rtimes_{\sigma,r} G
\egloz
(the $\varphi_i$ are given by \eqref{phii}) induces a KK-equivalence in $KK(\bigoplus_{i \in \cR} C^*_r(G_i), c_0(I) \rtimes_{\sigma,r} G)$. Moreover, we have
\bgl
\label{imp-eq}
  KK(\varphi_i) \otimes KK(\iota_{\eckl{i}}) \otimes j_r^G(\xf) = KK(\Phi_i)
\egl
where $\Phi_i$ is the homomorphism
\bgloz
  \Phi_i: C^*_r(G_i) \to D \rtimes_{\tau,r} G, \: \lambda_g \ma e_i U_g.
\egloz
\esubprop
\bproof
Since each of the $\varphi_i$ identifies $C^*_r(G_i)$ with a full corner of $C_{\eckl{i}}$, the homomorphism $\bigoplus_{i \in \cR} \varphi_i$ identifies the direct sum $\bigoplus_{i \in \cR} C^*_r(G_i)$ with a full corner of $\bigoplus_{i \in \cR} C_{\eckl{i}} = c_0(I) \rtimes_{\sigma,r} G$. This proves our first assertion. The second one is just \eqref{phix_i}.
\eproof

\section{From abstract to concrete}
\label{abstract-concrete}

Let us now go back to the situation of semigroup C*-algebras and summarize what we have obtained so far. We just have to apply our general results from the previous section to the case of reduced semigroup C*-algebras. We use the same notations as in \S~\ref{pre}. For a left Ore semigroup $P$ whose constructible right ideals are independent, set $P^{-1} \cdot (\cJ \setminus \gekl{\emptyset}) = \menge{q^{-1} \cdot X}{q \in P, \: \emptyset \neq X \in \cJ}$ and let $G$ be the enveloping group of $P$. The $G$-action on $P^{-1} \cdot (\cJ \setminus \gekl{\emptyset})$ via left multiplication (i.e. $g \cdot (q^{-1} \cdot X) = g \cdot q^{-1} \cdot X$) induces in a canonical way a $G$-action on $c_0(P^{-1} \cdot (\cJ \setminus \gekl{\emptyset}))$ by shifting indices. In the following, we will consider the conditions
\begin{itemize}
\item[(A${}_1$)] $P$ is a left Ore semigroup whose constructible right ideals are independent, and the enveloping group $G$ of $P$ satisfies the Baum-Connes conjecture with coefficients in the $G$-C*-algebras $c_0(P^{-1} \cdot (\cJ \setminus \gekl{\emptyset}))$ and $\DiP$,
\item[(A${}_2$)] Condition (A${}_1$) holds, and $c_0(P^{-1} \cdot (\cJ \setminus \gekl{\emptyset})) \rtimes_{\sigma,r} G$ and $\DiP \rtimes_{\taui,r} G$ satisfy the UCT or $G$ satisfies the strong Baum-Connes conjecture with coefficients in the $G$-C*-algebras $c_0(P^{-1} \cdot (\cJ \setminus \gekl{\emptyset}))$ and $\DiP$.
\end{itemize}

\btheo
\label{theosum}
If condition (A${}_1$) is satisfied, then the descent
\bgloz
  j_r^G(\xf) \in KK(c_0(P^{-1} \cdot (\cJ \setminus \gekl{\emptyset})) \rtimes_{\sigma,r} G, \DiP \rtimes_{\taui,r} G)
\egloz
of the element $\xf$ from Definition~\ref{def-xf} induces an isomorphism on K-theory.

If the stronger assumption (A${}_2$) is valid, then $j_r^G(\xf)$ is a KK-equivalence.
\etheo
\bproof
We have checked at the beginning of \S~\ref{concrete-abstract} that under the present assumptions, all the conditions in Theorem~\ref{mainthm} are satisfied. Hence the first part of the present theorem follows from Corollary~\ref{cor-mainthm1}, and the second part follows from Corollaries~\ref{cor-mainthm2} and \ref{cor-mainthm3}.
\eproof

Recall that the embedding $\iota: C^*_r(P) \to \DiP \rtimes_{\taui,r} G$, $V_p \ma  E_P U_p E_P$ induces a KK-equivalence in $KK(C^*_r(P), \DiP \rtimes_{\taui,r} G)$ by Corollary~\ref{CP-DG}. Also recall that for $X \in \cJ \setminus \gekl{\emptyset}$, we have introduced the homomorphism $\varphi_X: C^*_r(G_X) \to c_0(G \cdot X) \rtimes_{\sigma,r} G$, $\lambda_g \ma \1z_X U_g$ in \S~\ref{impthm}. Here $G_X = \menge{g \in G}{g \cdot X = X}$. Let $\iota_{\eckl{X}}$ be the embedding $c_0(G \cdot X) \rtimes_{\sigma,r} G \into c_0(P^{-1} \cdot \cJ \setminus \gekl{\emptyset}) \rtimes_{\sigma,r} G$.
\blemma
For every $X$ in $\cJ \setminus \gekl{\emptyset}$, there exists a homomorphism
\bgloz
  \Psi_X: C^*_r(G_X) \to C^*_r(P), \: \lambda_{q^{-1} p} \ma E_X V_q^* V_p E_X
\egloz
which satisfies
\bgl
\label{Psi}
  KK(\Psi_X) = KK(\varphi_X) \otimes KK(\iota_{\eckl{X}}) \otimes j_r^G(\xf) \otimes KK(\iota)^{-1}.
\egl
\elemma
\bproof
Let $X$ be an element of $\cJ \setminus \gekl{\emptyset}$. Recall that $\Phi_X$ is the homomorphism $C^*_r(G_X) \to \DiP \rtimes_{\taui,r} G, \lambda_g \ma E_X U_g = E_X U_g E_X$. It is clear that $\img(\Phi_X) \subseteq \img(\iota)$, so that we can define $\Psi_X \defeq \iota^{-1} \circ (\Phi_X \vert^{\img(\iota)})$. This homomorphism has the desired properties. Equation~\eqref{Psi} follows from $\iota \circ \Psi_X = \Phi_X$ (by construction) and \eqref{imp-eq}.
\eproof

Now let $\cX$ be a complete system of representatives for $G \backslash \rukl{P^{-1} \cdot (\cJ \setminus \gekl{\emptyset})}$ such that $\cX \subseteq \cJ \setminus \gekl{\emptyset}$. The homomorphisms $\gekl{\Psi_X}_{X \in \cX}$ from the previous lemma give rise to the Kasparov $(\bigoplus_{X \in \cX} C^*_r(G_X),C^*_r(P))$-module $(\ell^2(\cX, C^*_r(P)), \bigoplus_{X \in \cX} \Psi_X, 0)$ with the homomorphism
\bgloz
  \bigoplus_{X \in \cX} \Psi_X: 
  \bigoplus_{X \in \cX} C^*_r(G_X) \to c_0(\cX, C^*_r(P)) \subseteq \cK(\ell^2(\cX, C^*_r(P))) \subseteq \cL(\ell^2(\cX, C^*_r(P))).
\egloz
Here $c_0(\cX, C^*_r(P))$ acts as diagonal multiplication operators on the Hilbert $C^*_r(P)$-module $\ell^2(\cX, C^*_r(P))$. Let $\Psi$ be the KK-element in $KK(\bigoplus_{X \in \cX} C^*_r(G_X),C^*_r(P))$ represented by the Kasparov module $(\ell^2(\cX, C^*_r(P)), \bigoplus_{X \in \cX} \Psi_X, 0)$. Let $\iota_{G_X}$ be the inclusion $C^*_r(G_X) \into \bigoplus_{X \in \cX} C^*_r(G_X)$. By construction, 
\bgl
\label{PsiPsi}
  KK(\iota_{G_X}) \otimes \Psi = KK(\Psi_X).
\egl

\btheo
\label{PsiKK}
If condition (A${}_1$) is valid, then the KK-element $\Psi$ from above induces an isomorphism on K-theory.

If the stronger assumption (A${}_2$) holds, then $\Psi$ is a KK-equivalence.
\etheo
\bproof
By Corollary~\ref{CP-DG}, $KK(\iota)$ is a KK-equivalence. By the first part of Proposition~\ref{imp}, $KK(\bigoplus_{X \in \cX} \varphi_X)$ is a KK-equivalence. And going through the identification
\bgloz
  KK(\bigoplus_{X \in \cX} C^*_r(G_X),C^*_r(P)) \cong \prod_{X \in \cX} KK( C^*_r(G_X),C^*_r(P))
\egloz
from \cite{Bla}, Theorem~19.7.1, it follows from equation~\eqref{Psi} from the previous lemma and \eqref{PsiPsi} that
\bgloz
  \Psi = KK(\bigoplus_{X \in \cX} \varphi_X) \otimes j_r^G(\xf) \otimes KK(\iota)^{-1}.
\egloz

Therefore, the first part of the present theorem follows from the first part of Theorem~\ref{theosum}, and the second part follows from the second part of the same theorem.
\eproof

\bcor
\label{PsiKKcor}
If condition (A${}_1$) is satisfied, then the homomorphism
\bgloz
  \sum_{X \in \cX} K_*(\Psi_X): \bigoplus_{X \in \cX} K_*(C^*_r(G_X)) \to K_*(C^*_r(P))
\egloz
is an isomorphism. And under the stronger assumption (A${}_2$), the homomorphism
\bgloz
  \prod_{X \in \cX} K^*(\Psi_X): K^*(C^*_r(P)) \to \prod_{X \in \cX} K^*(C^*_r(G_X))
\egloz
is an isomorphism. Here $K_*$ is $K_0 \oplus K_1$ and $K^*$ is $K^0 \oplus K^1$ viewed as $\Zz / 2 \Zz$-graded abelian groups.
\ecor
\bproof
By \eqref{PsiPsi}, these homomorphisms are just the compositions of taking the Kasparov product with $\Psi$ and the canonical isomorphisms
\bgloz
  \bigoplus_{X \in \cX} K_*(C^*_r(G_X)) \cong K_*(\bigoplus_{X \in \cX} C^*_r(G_X)) 
  \text{ and } 
  K^*(\bigoplus_{X \in \cX} C^*_r(G_X)) \cong \prod_{X \in \cX} K^*(C^*_r(G_X)).
\egloz
\eproof

Our last goal in this section is to show that whenever $P$ is a left Ore semigroup whose constructible right ideals are independent, there exists a canonical ring structure on the K-homology of $C^*_r(P)$ and that the isomorphism $\prod_{X \in \cX} K^*(\Psi_X)$ from the last corollary is a ring isomorphism.

\blemma
\label{diagex}
Let $P$ be a left Ore semigroup whose constructible right ideals are independent. Then there exists a homomorphism
\bgloz
  \Delta_P: C^*_r(P) \to C^*_r(P) \otimes_{\min} C^*_r(P) \text{ determined by } V_p \ma V_p \otimes V_p.
\egloz
\elemma

Note that we always have such a homomorphism in the case where the left regular representation $C^*(P) \to C^*_r(P)$ is an isomorphism because an analogous homomorphism always exists on the full semigroup C*-algebra (see the proof of Proposition~2.24 in \cite{Li2}).

\bproof
Since the constructible right ideals of $P$ are independent, there exists a homomorphism $D_r(P) \to D_r(P) \otimes_{\min} D_r(P)$ sending $E_X$ to $E_X \otimes E_X$ for all $X \in \cJ$. This can be seen as follows: By \cite{Li2}, Corollary~2.26, the restriction of the left regular representation to the commutative sub-C*-algebra $D(P)$ of the full semigroup C*-algebra $C^*(P)$ yields an isomorphism $D(P) \cong D_r(P)$ if (and only if) the constructible right ideals of $P$ are independent. But we can always construct a homomorphism $D(P) \to D(P) \otimes_{\min} D(P)$, $e_X \ma e_X \otimes e_X$ by restricting the homomorphism $C^*(P) \to C^*(P) \otimes_{\min} C^*(P)$, $v_p \ma v_p \otimes v_p$ to $D(P)$ (as observed above, such a homomorphism always exists; see also the proof of Proposition~2.24 in \cite{Li2}).

The homomorphism $D_r(P) \to D_r(P) \otimes_{\min} D_r(P)$, $E_X \ma E_X \otimes E_X$ is obviously equivariant with respect to the $P$-actions $\tau$ and $\tau \otimes \tau$. By definition of $\DiP$ (see the beginning of \S~\ref{dilations}), we obtain a homomorphism $\DiP \to \DiP \otimes_{\min} \DiP$, $E_Y \ma E_Y \otimes E_Y$ (for $Y \in P^{-1} \cdot \cJ$). This homomorphism is again obviously $G$-equivariant with respect to the actions $\taui$ and $\taui \otimes \taui$. Therefore, applying Lemma~\ref{functor} to this homomorphism and the diagonal embedding $G \into G \times G$, we obtain the homomorphism
\bglnoz
  \DiP \rtimes_{\taui,r} G &\to& (\DiP \otimes_{\min} \DiP) \rtimes_{\taui \otimes \taui,r} (G \times G) \\
  E_Y U_g &\ma& (E_Y \otimes E_Y) U_{(g,g)}.
\eglnoz
Composing this map with the canonical identification
\bglnoz
  && (\DiP \otimes_{\min} \DiP) \rtimes_{\taui \otimes \taui,r} (G \times G) \\
  &\cong& (\DiP \rtimes_{\taui,r} G) \otimes_{\min} (\DiP \rtimes_{\taui,r} G), \\
  && (E_{Y_1} \otimes E_{Y_2}) U_{(g_1,g_2)} \ma E_{Y_1} U_{g_1} \otimes E_{Y_2} U_{g_2}
\eglnoz
we obtain the homomorphism
\bglnoz
  \DiP \rtimes_{\taui,r} G &\to& (\DiP \rtimes_{\taui,r} G) \otimes_{\min} (\DiP \rtimes_{\taui,r} G)\\
  E_Y U_g &\ma& E_Y U_g \otimes E_Y U_g.
\eglnoz
Since this map sends $E_P U_p E_P$ to $E_P U_p E_P \otimes E_P U_p E_P$, we just have to restrict this homomorphism to $E_P \rukl{\DiP \rtimes_{\taui,r} G} E_P$ and to use the identification 
\bgloz
  C^*_r(P) \cong E_P \rukl{\DiP \rtimes_{\taui,r} G} E_P, \: V_p \ma E_P U_p E_P
\egloz
from Lemma~\ref{P-EGE} to obtain our desired homomorphism $\Delta_P$.
\eproof

Now, whenever there exists such a diagonal homomorphism
\bgloz
  \Delta_P: C^*_r(P) \to C^*_r(P) \otimes_{\min} C^*_r(P), \: V_p \ma V_p \otimes V_p,
\egloz
we obtain a canonical graded ring structure on $K^*(C^*_r(P))$ in analogy to the group case. Multiplication in this ring structure is given by the following composition:
\bglnoz
  K^i(C^*_r(P)) \times K^j(C^*_r(P)) &\cong& KK^i(C^*_r(P),\Cz) \times KK^j(C^*_r(P),\Cz) \\
  &\overset{\otimes}{\lori}& KK^{i+j}(C^*_r(P) \otimes_{\min} C^*_r(P),\Cz) \\
  &\overset{KK(\Delta_P) \otimes \sqcup}{\lori}& KK^{i+j}(C^*_r(P),\Cz) \cong K^{i+j}(C^*_r(P)).
\eglnoz
And on $\prod_{X \in \cX} K^*(C^*_r(G_X))$, there is a canonical ring structure given by the canonical ring structure on each of the K-homology groups $K^*(C^*_r(G_X))$ (it is constructed in the same way as it was explained for $K^*(C^*_r(P))$). Our last observation in this section is that the isomorphism on K-homology from the last corollary is compatible with these ring structures.
\btheo
\label{ringst}
If condition (A${}_2$) is satisfied, then the homomorphism
\bgloz
  \prod_{X \in \cX} K^*(\Psi_X): K^*(C^*_r(P)) \to \prod_{X \in \cX} K^*(C^*_r(G_X))
\egloz
is a ring isomorphism.
\etheo
\bproof
In view of the last corollary, all we have to prove is that $\prod_{X \in \cX} K^*(\Psi_X)$ is multiplicative. Let us check this for $K^0$, the remaining cases are similar. Let $\Delta_{G_X}$ be the diagonal homomorphism $C^*_r(G_X) \to C^*_r(G_X) \otimes_{\min} C^*_r(G_X), \lambda_g \ma \lambda_g \otimes \lambda_g$. Using the natural identification $K^0(\cdot) \cong KK(\cdot,\Cz)$ and the definition of the multiplicative structures, our assertion amounts to saying that for all $X \in \cX$ and for all $\yf$, $\zf$ in $KK(C^*_r(P),\Cz)$, we have
\bgloz
  KK(\Psi_X) \otimes (KK(\Delta_P) \otimes (\yf \otimes \zf)) = KK(\Delta_{G_X}) \otimes \rukl{(KK(\Psi_X) \otimes \yf) \otimes (KK(\Psi_X) \otimes \zf)}.
\egloz
But by \cite{Kas}, Theorem~2.14~8), we obtain
\bgl
\label{commrel}
  \yf \otimes KK(\Psi_X) = KK(\Psi_X) \otimes \yf.
\egl
Moreover, it is immediate that
\bgl
\label{diagPsi}
  \Delta_P \circ \Psi_X = (\Psi_X \otimes_{\min} \Psi_X) \circ \Delta_{G_X}.
\egl
Thus
\bglnoz
  && KK(\Delta_{G_X}) \otimes \rukl{(KK(\Psi_X) \otimes \yf) \otimes (KK(\Psi_X) \otimes \zf)} \\
  &\overset{\eqref{commrel}}{=}& KK(\Delta_{G_X}) \otimes (KK(\Psi_X) \otimes KK(\Psi_X)) \otimes (\yf \otimes \zf) \\
  &=& KK(\Delta_{G_X}) \otimes KK(\Psi_X \otimes_{\min} \Psi_X) \otimes (\yf \otimes \zf) \\
  &\overset{\eqref{diagPsi}}{=}& KK(\Psi_X) \otimes (KK(\Delta_P) \otimes (\yf \otimes \zf)).
\eglnoz
\eproof

\section{Semigroups attached to Dedekind domains}
\label{Dede}

In this section, we apply our general K-theoretic results from \S~\ref{abstract-concrete} to specific semigroups attached to Dedekind domains. Let $R$ be a Dedekind domain. This means that $R$ is a noetherian, integrally closed integral domain with the property that every non-zero prime ideal is a maximal ideal (compare \cite{Neu}, Chapter~I, Definition~(3.2)). By an integral domain, we mean a commutative ring without zero divisors.

We would like to treat the multiplicative semigroup $R\reg = R \setminus \gekl{0}$, the semigroup of principal ideals of $R$ and the $ax+b$-semigroup $R \rtimes R\reg$. The semidirect product $R \rtimes R\reg$ is taken with respect to the multiplicative action of the multiplicative semigroup $R\reg$ on the additive group $R$.

Examples of Dedekind domains are given by rings of integers in number fields or function fields. These rings and the corresponding semigroups have actually been our motivating examples.

Since it will be important later on, let us briefly recall the definition of the class group of $R$. Let $Q(R)$ be the quotient field of $R$.
\bdefin
A fractional ideal of $Q(R)$ (or $R$) is a non-zero, finitely generated sub-$R$-module of $Q(R)$.

A principal fractional ideal of $Q(R)$ (or $R$) is a fractional ideal of the form $a \cdot R$ for some $a \in Q(R)\reg = Q(R) \setminus \gekl{0}$.
\edefin
As explained in \cite{Neu}, Chapter~I, \S~3, the set of fractional ideals of $Q(R)$ form an abelian group under multiplication. Furthermore, the subset of principal fractional ideals of $Q(R)$ is multiplicatively closed, hence it forms a subgroup.
\bdefin
The ideal class group (or simply class group) $Cl_{Q(R)}$ of $Q(R)$ is the quotient of the group of fractional ideals by the subgroup of principal fractional ideals of $Q(R)$.
\edefin
\bremark
\label{Cl}
It follows directly from the definition that we can equivalently describe $Cl_{Q(R)}$ (at least as a set) as follows: The multiplicative group $Q(R)\reg = Q(R) \setminus \gekl{0}$ acts on the set of fractional ideals of $Q(R)$ by multiplication, and $Cl_{Q(R)}$ is given by the set of orbits of this action.
\eremark

\subsection{Multiplicative semigroups}

We first consider the multiplicative semigroup $R\reg$. Let $R^*$ be the group of units in $R$, or in other words, $R^*$ is the subgroup of invertible elements of $R\reg$.

Our goal is to prove as an application of our general K-theoretic results from \S~\ref{abstract-concrete}:
\bsubtheo
$C^*_r(R\reg)$ and $\bigoplus_{\gamma \in Cl_{Q(R)}} C^*_r(R^*)$ are KK-equivalent.

Furthermore, choose for every $\gamma \in Cl_{Q(R)}$ an ideal $I_{\gamma}$ of $R$ which represents $\gamma$. Then there is a homomorphism $\Psi_{I_{\gamma}}: C^*_r(R^*) \to C^*(R\reg)$ determined by $\Psi_{I_{\gamma}}(V_a) = E_{I_{\gamma}} V_a$. These homomorphisms give rise to isomorphisms
\bgloz
  \sum_{\gamma \in Cl_{Q(R)}} (\Psi_{I_{\gamma}})_*: \bigoplus_{\gamma \in Cl_{Q(R)}} K_*(C^*_r(R^*)) \to K_*(C^*_r(R\reg))
\egloz
and
\bgloz
  \prod_{\gamma \in Cl_{Q(R)}} (\Psi_{I_{\gamma}})^*: K^*(C^*_r(R\reg)) \to \prod_{\gamma \in Cl_{Q(R)}} K^*(C^*_r(R^*)).
\egloz
The last isomorphism $\prod_{\gamma \in Cl_{Q(R)}} (\Psi_{I_{\gamma}})^*$ on K-homology is a ring isomorphism.
\esubtheo
\bproof
We just have to check the assumptions in Theorem~\ref{PsiKK}. First of all, $R\reg$ is a left Ore semigroup because it is cancellative and abelian. Moreover, the constructible right ideals of $R\reg$ are independent. This can be proven analogously to Lemma~2.29 in \cite{Li2}. The enveloping group of $R\reg$ is given by $Q(R)\reg$. Since $Q(R)\reg$ is abelian, it is amenable, hence it satisfies the strong Baum-Connes conjecture for all coefficients. Therefore the conditions in the second part of Theorem~\ref{PsiKK} are satisfied. For the semigroup $R\reg$, $\cJ \setminus \gekl{\emptyset}$ is given by all non-zero ideals of $R$. This can be proven analogously to the case of the $ax+b$-semigroup over $R$ which is explained in the second half of \S~2.4 of \cite{Li2}. Therefore, $(R\reg)^{-1} \cdot (\cJ \setminus \gekl{\emptyset})$ is the set of fractional ideals of $Q(R)$, and the set of orbits $Q(R)\reg \backslash (P^{-1} \cdot (\cJ \setminus \gekl{\emptyset}))$ coincides with $Cl_{Q(R)}$ by Remark~\ref{Cl}. And finally, for a non-zero ideal $I$ of $R$, the stabilizer $Q(R)\reg_I = \menge{a \in Q(R)\reg}{a \cdot I = I}$ is given by $R^*$. The first part of our theorem now follows from the second part of Theorem~\ref{PsiKK} and from the second part of Corollary~\ref{PsiKKcor}. That $\prod_{\gamma \in Cl_{Q(R)}} (\Psi_{I_{\gamma}})^*$ is a ring isomorphism follows from Theorem~\ref{ringst}.
\eproof

\bsubremark
Let $L$ be an ideal in $R$. We define an (non-unital) endomorphism $\alpha_L$ of $C^*_r(R^\times)$ by $V_p \mapsto V_p E_L$, $E_I\mapsto E_{LI}$. Then one has $\alpha_J\alpha_L=\alpha_{JL}$ and $\alpha_L$ is inner if $L$ is a principal ideal.

As a consequence we obtain an action of the class group $Cl_{Q(R)}$ on the K-theory and K-homology of $C^*_r(R^\times)$ (in fact this defines a multiplicative map $Cl_{Q(R)} \to KK(C^*_r(R^\times),C^*_r(R^\times))$). It is clear that this action of $Cl_{Q(R)}$ corresponds under $\sum_{\gamma \in Cl_{Q(R)}} (\Psi_{I_{\gamma}})_*$ to the obvious action of $Cl_{Q(R)}$ on $\bigoplus_{\gamma \in Cl_{Q(R)}} K_*(C^*_r(R^*))$, and similarly on K-homology.
\esubremark

We also discuss the multiplicative semigroup of principal ideals over a Dedekind domain $R$. It is clear that this semigroup can be identified with $R\reg / R^*$. Note that the family of ideals $\cJ_{R\reg / R^*}$ for this semigroup can be identified with the corresponding family $\cJ_{R\reg}$ for the multiplicative semigroup of the ring $R$ via $\cJ_{R\reg / R^*} \ni X / R^* \leftrightarrow X \in \cJ_{R\reg}$, where $X / R^*$ is the image of $X$ in $R\reg / R^*$ under the canonical projection $R\reg \onto R\reg / R^*$. With this observation, we can in complete analogy to the case of $R\reg$ apply Theorem~\ref{PsiKK}, Corollary~\ref{PsiKKcor} and Theorem~\ref{ringst} and deduce
\bsubtheo
$C^*_r(R\reg / R^*)$ and $\bigoplus_{\gamma \in Cl_{Q(R)}} \Cz$ are KK-equivalent.

Furthermore, choose for every $\gamma \in Cl_{Q(R)}$ an ideal $I_{\gamma}$ of $R$ which represents $\gamma$. Then the canonical homomorphism $\Psi_{I_{\gamma}}: \Cz \to C^*(R\reg / R^*)$ determined by $\Psi_{I_{\gamma}}(1) = E_{I_{\gamma} / R^*}$ give rise to isomorphisms
\bgloz
  \sum_{\gamma \in Cl_{Q(R)}} (\Psi_{I_{\gamma}})_*: \bigoplus_{\gamma \in Cl_{Q(R)}} \Zz \to K_*(C^*_r(R\reg))
\egloz
and
\bgloz
  \prod_{\gamma \in Cl_{Q(R)}} (\Psi_{I_{\gamma}})^*: K^*(C^*_r(R\reg)) \to \prod_{\gamma \in Cl_{Q(R)}} \Zz.
\egloz
The last isomorphism $\prod_{\gamma \in Cl_{Q(R)}} (\Psi_{I_{\gamma}})^*$ is a ring isomorphism, where we take the canonical ring structure on $\Zz$.
\esubtheo

We also remark that we obtain an analogous action of the class group $Cl_{Q(R)}$ on the K-theory and K-homology of $C^*_r(R^\times / R^*)$ as in the previous remark.

\subsection{$ax+b$-semigroups}
\label{ax+b}

Let us now treat the case of the $ax+b$-semigroup $R \rtimes R\reg$ over $R$. First, we apply our general results from \S~\ref{abstract-concrete} in order to compute K-theory, and secondly, we show that the corresponding semigroup C*-algebras are purely infinite.

Again, let $R^*$ be the group of units in $R$ and choose for every $\gamma \in Cl_{Q(R)}$ an ideal $I_{\gamma}$ of $R$ which represents $\gamma$.

Applying our general K-theoretic results from \S~\ref{abstract-concrete}, we obtain
\bsubtheo
\label{Kax+b}
The C*-algebras $C^*_r(R \rtimes R\reg)$ and $\bigoplus_{\gamma \in Cl_{Q(R)}} C^*_r(I_{\gamma} \rtimes R^*)$ are KK-equivalent. Here we form the semidirect product $I_{\gamma} \rtimes R^*$ with respect to the multiplicative action of $R^*$ on the additive group $I_{\gamma}$.

Moreover, for every $\gamma \in Cl_{Q(R)}$ there is a homomorphism $\Psi_{I_{\gamma}}: C^*_r(I_{\gamma} \rtimes R^*) \to C^*(R \rtimes R\reg)$ determined by $\Psi_{I_{\gamma}}(V_{(b,a)}) = E_{I_{\gamma} \times I_{\gamma}\reg} V_{(b,a)}$. These homomorphisms give rise to isomorphisms
\bgloz
  \sum_{\gamma \in Cl_{Q(R)}} (\Psi_{I_{\gamma}})_*: \bigoplus_{\gamma \in Cl_{Q(R)}} K_*(C^*_r(I_{\gamma} \rtimes R^*)) \to K_*(C^*_r(R \rtimes R\reg))
\egloz
and
\bgloz
  \prod_{\gamma \in Cl_{Q(R)}} (\Psi_{I_{\gamma}})^*: K^*(C^*_r(R \rtimes R\reg)) \to \prod_{\gamma \in Cl_{Q(R)}} K^*(C^*_r(I_{\gamma} \rtimes R^*)).
\egloz
The last isomorphism $\prod_{\gamma \in Cl_{Q(R)}} (\Psi_{I_{\gamma}})^*$ on K-homology is a ring isomorphism.
\esubtheo
\bproof
Again, we just have to check the assumptions in Theorem~\ref{PsiKK}. First of all, $R \rtimes R\reg$ is a left Ore semigroup by \cite{Li1}, \S~5.1. And the constructible right ideals of $R \rtimes R\reg$ are independent by Lemma~2.29 in \cite{Li2}. The enveloping group of $R \rtimes R\reg$ is given by the $ax+b$-group $Q(R) \rtimes Q(R)\reg$ over $Q(R)$. Since $Q(R) \rtimes Q(R)\reg$ is solvable, it is amenable, hence it satisfies the strong Baum-Connes conjecture for all coefficients. Therefore the conditions in the second part of Theorem~\ref{PsiKK} are fulfilled. For the semigroup $R \rtimes R\reg$, $\cJ \setminus \gekl{\emptyset}$ is given by $\menge{(r+I) \times I\reg}{r \in R, (0) \neq I \triangleleft R}$. This is explained in the second half of \S~2.4 of \cite{Li2}. Therefore,
\bgloz
  (R \rtimes R\reg)^{-1} \cdot (\cJ \setminus \gekl{\emptyset}) = \menge{(a^{-1}b + a^{-1}I) \times (a^{-1} I)\reg}{(b,a) \in R \rtimes R\reg, (0) \neq I \triangleleft R},
\egloz
and we see using Remark~\ref{Cl} that
\bgloz
  Cl_{Q(R)} \to (Q(R) \rtimes Q(R)\reg) \backslash ((R \rtimes R\reg)^{-1} \cdot (\cJ \setminus \gekl{\emptyset})), \: J \ma \eckl{J \times J\reg}
\egloz
is a bijection. And finally, for a non-zero ideal $I$ of $R$, the stabilizer $(Q(R) \rtimes Q(R)\reg)_{I \times I\reg} = \menge{(b,a) \in Q(R) \rtimes Q(R)\reg}{b + a \cdot I = I}$ is given by $I \rtimes R^* \subseteq R \rtimes R\reg$. The first part of our theorem now follows from the second part of Theorem~\ref{PsiKK} and the second part of Corollary~\ref{PsiKKcor}, and Theorem~\ref{ringst} implies that $\prod_{\gamma \in Cl_{Q(R)}} (\Psi_{I_{\gamma}})^*$ is a ring isomorphism.
\eproof

Finally, let us study the inner structure of semigroup C*-algebras of $ax+b$-semigroups over Dedekind domains. We start with two definitions:

\bsubdefin
A C*-algebra $A$ is purely infinite if $A$ has no non-zero abelian quotients and for every pair of positive elements $a$ and $b$ in $A$ with $b \in \overline{AaA}$, there exists a sequence $(x_n)_n$ in $A$ such that $\lim_{n \to \infty} x_n^* a x_n = b$.
\esubdefin
The reader may consult \cite{Ror}, \cite{Pas-Ror} or \cite{Kir-Ror1} for more details.

\bsubdefin
A C*-algebra has the ideal property if projections separate ideals.
\esubdefin
Further explanations can be found in \cite{Pas-Ror}.

Our final goal is to prove
\bsubtheo
\label{pi}
For every Dedekind domain which has infinitely many pairwise distinct prime ideals, the semigroup C*-algebra $C^*_r(R \rtimes R\reg)$ is purely infinite and has the ideal property.
\esubtheo

For us, the following result of C. Pasnicu and M. R{\o}rdam (see \cite{Pas-Ror}, Proposition~2.11) is important:

A C*-algebra is purely infinite and has the ideal property if and only if every non-zero hereditary sub-C*-algebra in any quotient contains an infinite projection.

Actually, we will only need the implication \an{$\Larr$}. Our goal is to prove that for every ideal $\fI$ of $C^*_r(R \rtimes R\reg)$, every non-zero hereditary sub-C*-algebra of $C^*_r(R \rtimes R\reg) / \fI$ contains an infinite projection.

Let us start with a general observation. Let $D$ be a unital C*-algebra with an action $\alpha$ of a semigroup $P$ by injective endomorphisms. Form the semigroup crossed product $D \rte_{\alpha} P$ in the sense of \cite{La} or \cite{Li1}, \S~A1. Recall that $D \rte_{\alpha} P$ is a unital C*-algebra which comes by definition with a unital homomorphism $i_D: D \to D \rte_{\alpha} P$ and a semigroup homomorphism $v: P \to \Isom(D \rte_{\alpha} P)$, $p \ma v_p$ such that $v_p i_D(d) v_p^* = i_D(\alpha_p(d))$ for all $p \in P$ and $d \in D$. The triple $(D \rte_{\alpha} P,i_D,v)$ satisfies the universal property that given a unital C*-algebra $T$, a unital homomorphism $j_D: D \to T$ and a semigroup homomorphism $w: P \to \Isom(T)$ such that $w_p j_D(d) w_p^* = j_D(\alpha_p(d))$ for all $p \in P$ and $d \in D$, there is a unique homomorphism $j_D \rte w: D \rte_{\alpha} P \to T$ satisfying $(j_D \rte w) \circ i_D = j_D$ and $(j_D \rte w) \circ v = w$.

\bsublemma
\label{pi1}
Assume that $P$ is a left Ore semigroup and that the enveloping group $G$ of $P$ is amenable. Moreover, assume that
\bgl
\label{sat}
  v_p^* i_D(d) v_p \in i_D(D) \fa d \in D.
\egl
Then there exists a faithful conditional expectation $E: D \rte_{\alpha} P \to D$ which is uniquely determined by
\bgl
\label{E-det}
E(v_q^* i_D(d) v_p) = \delta_{q,p} i_D^{-1}(v_p^* i_D(d) v_p).
\egl
\esublemma
\bproof
In the situation of the lemma, we have
\bgloz
  D \rte_{\alpha} P = \clspan{\menge{v_q^* i_D(d) v_p}{p,q \in P; d \in D}}
\egloz
by Remark~1.3.1 in \cite{La}. This explains why \eqref{E-det} completely determines $E$.

To prove existence of $E$, let $(D_{\infty},G,\alpha\iu)$ be the minimal automorphic dilation in the sense of Definition~2.1.2 of \cite{La}. By Theorem~2.2.1 in \cite{La}, we have canonical embeddings $i: D \into D_{\infty}$, $i^{(\rtimes)}: D \rte_{\alpha} P \into D_{\infty} \rtimes_{\alpha\iu} G$ and $D_{\infty} \overset{\subseteq}{\into} D_{\infty} \rtimes_{\alpha\iu} G$ such that the diagram
\bgloz
  \begin{CD}
  D \rte_{\alpha} P @> i^{(\rtimes)} >> D_{\infty} \rtimes_{\alpha\iu} G \\
  @AA i_D A @AA \subseteq A \\
  D @> i >> D_{\infty}
  \end{CD}
\egloz
commutes. It follows that $i_D$ is injective. Moreover, Theorem~2.2.1 in \cite{La} tells us that $\img(i^{(\rtimes)}) = i(1_D) \rukl{D_{\infty} \rtimes_{\alpha\iu} G} i(1_D)$.

Now, as $G$ is amenable, $D_{\infty} \rtimes_{\alpha\iu} G \cong D_{\infty} \rtimes_{\alpha\iu,r} G$. And since $G$ is also discrete, Lemma~\ref{fce} implies that there is a faithful conditional expectation $E_{\infty}: D_{\infty} \rtimes_{\alpha\iu} G \to D_{\infty}$ determined by
\bgloz
  E_{\infty} (d_{\infty} u_g) = \delta_{g,e} d_{\infty}.
\egloz
Here $u_g$ are the canonical unitaries in the multiplier algebra of $D_{\infty} \rtimes_{\alpha\iu} G$ which implement $\alpha\iu$. As $E_{\infty}(i(1_D)) = i(1_D)$, the composition
\bgloz
  D \rte_{\alpha} P \overset{i^{(\rtimes)}}{\lori} D_{\infty} \rtimes_{\alpha\iu} G \overset{E_{\infty}}{\lori} D_{\infty}
\egloz
has image in $i(1_D) \rukl{D_{\infty} \rtimes_{\alpha\iu} G} i(1_D) = \img(i^{(\rtimes)})$, so that we can form
\bgloz
  E' \defeq (i^{(\rtimes)})^{-1} \circ (E_{\infty} \circ i^{(\rtimes)}) \vert^{\img(i^{(\rtimes)})}.
\egloz
$E'$ is a faithful conditional expectation determined by
\bgl
\label{E'}
  E'(v_q^* i_D(d) v_p) = \delta_{q,p} v_p^* i_D(d) v_p.
\egl
By our assumption \eqref{sat}, $v_p^* i_D(d) v_p$ lies in $i_D(D)$ for all $p \in P$ and $d \in D$. Thus $\img(E') = i_D(D)$, and since $i_D$ is injective, we may set
\bgloz
  E \defeq i_D^{-1} \circ (E' \vert^{i_D(D)}).
\egloz
This is the desired faithful conditional expectation. It satisfies \eqref{E-det} because of \eqref{E'}.
\eproof

\bsublemma
\label{pi2}
Let us consider the same situation as in the previous lemma. Let $\fI$ be an ideal of $D \rte_{\alpha} P$. Then
\bgloz
  \fI_D \defeq i_D^{-1}(i_D(D) \cap \fI)
\egloz
is an ideal of $D$ such that for every $p \in P$, the endomorphism
\bgloz
  \dot{\alpha}_p: D / \fI_D \to D / \fI_D, d + \fI_D \ma \alpha_p(d) + \fI_D
\egloz
is well-defined and injective. Let us denote the corresponding $P$-action on $D / \fI_D$ by $\dot{\alpha}$ and the associated semigroup crossed product by $((D / \fI_D) \rte_{\dot{\alpha}} P,i_{D / \fI_D},\dot{v})$. Let $\pi$ be the canonical projection $D \onto D / \fI_D$. By universal property of $(D \rte_{\alpha} P,i_D,v)$, there exists a homomorphism $\pi \rte P: D \rte_{\alpha} P \to (D / \fI_D) \rte_{\dot{\alpha}} P$ determined by
\bgloz
  (\pi \rte P)(v_q^* i_D(d) v_p) = \dot{v}_q^* i_{D / \fI_D}(d + \fI_D) \dot{v}_p.
\egloz
This homomorphism $\pi \rte P$ induces an isomorphism
\bgloz
  (\pi \rte P)\dot{}: D \rte_{\alpha} P / \spkl{i_D(\fI_D)} \overset{\cong}{\lori} (D / \fI_D) \rte_{\dot{\alpha}} P
\egloz
determined by
\bgl
\label{piP}
  (\pi \rte P)\dot{} \, (v_q^* i_D(d) v_p + \spkl{i_D(\fI_D)}) = \dot{v}_q^* i_{D / \fI_D}(d + \fI_D) \dot{v}_p.
\egl
Here $\spkl{i_D(\fI_D)}$ is the ideal of $D \rte_{\alpha} P$ generated by $i_D(\fI_D)$.
\esublemma
\bproof
If $d$ lies in $\fI_D$, then $i_D(\alpha_p(d)) = v_p(i_D(d))v_p^*$ lies in $\fI$ as $i_D(d)$ lies in $\fI$. At the same time, $v_p(i_D(d))v_p^* = i_D(\alpha_p(d))$ lies in $i_D(D)$. Thus $i_D(\alpha_p(d))$ lies in $i_D(D) \cap \fI$, and hence $\alpha_p(d)$ lies in $i_D^{-1}(i_D(D) \cap \fI) = \fI_D$. Therefore $\dot{\alpha}_p$ is well-defined. To see injectivity of $\alpha_p$, observe that for $d \in D$, $\alpha_p(d) \in \fI_D$ implies
\bgloz
  d = i_D^{-1} i_D(d) = i_D^{-1}(v_p^* v_p i_D(d) v_p^* v_p) = i_D^{-1}(v_p^* i_D(\alpha_p(d)) v_p).
\egloz
Now $v_p^* i_D(\alpha_p(d)) v_p$ lies in $\fI$ as $i_D(\alpha_p(d))$ lies in $i_D(\fI_D) \subseteq \fI$, and $v_p^* i_D(\alpha_p(d)) v_p$ lies in $i_D(D)$ by \eqref{sat}. Hence $v_p^* i_D(\alpha_p(d)) v_p$ lies in $i_D(D) \cap \fI$, and thus $d = i_D^{-1}(v_p^* i_D(\alpha_p(d)) v_p)$ lies in $i_D^{-1}(i_D(D) \cap \fI) = \fI_D$. So far, we have proven the first part of the lemma.

Finally, the homomorphism $\pi \rte P: D \rte_{\alpha} P \to (D / \fI_D) \rte_{\dot{\alpha}} P$ determined by
\bgloz
  (\pi \rte P)(v_q^* i_D(d) v_p) = \dot{v}_q^* i_{D / \fI_D}(d + \fI_D) \dot{v}_p
\egloz
vanishes on $i_D(\fI_D)$, hence on $\spkl{i_D(\fI_D)}$. Therefore $\pi \rte P$ factorizes through the quotient $D \rte_{\alpha} P / \spkl{i_D(\fI_D)}$. This gives rise to the desired homomorphism $(\pi \rte P)\dot{}$. To prove that $(\pi \rte P)\dot{}$ is an isomorphism, we use the universal property of $((D / \fI_D) \rte_{\dot{\alpha}} P, i_{D / \fI_D}, \dot{v})$ to construct an inverse. The composition
\bgloz
  D \overset{i_D}{\lori} D \rte_{\alpha} P \overset{\pi^{(\rtimes)}}{\onto} D \rte_{\alpha} P / \spkl{i_D(\fI_D)}
\egloz
($\pi^{(\rtimes)}$ is the canonical projection) obviously vanishes on $\fI_D$, so that we obtain a homomorphism $(i_D)\dot{}: D / \fI_D \to D \rte_{\alpha} P / \spkl{i_D(\fI_D)}$. It is straightforward to see that $(i_D)\dot{}$ and $P \ni p \ma \pi^{(\rtimes)}(v_p) \in \Isom(D \rte_{\alpha} P / \spkl{i_D(\fI_D)})$ satisfy the covariance relation $\Ad(\pi^{(\rtimes)}(v_p)) \circ (i_D)\dot{} = (i_D)\dot{} \circ \dot{\alpha}_p$. Therefore, by universal property of $((D / \fI_D) \rte_{\dot{\alpha}} P, i_{D / \fI_D}, \dot{v})$, there exists a homomorphism
\bgloz
  (i_D)\dot{} \rte P: (D / \fI_D) \rte_{\dot{\alpha}} P \to D \rte_{\alpha} P / \spkl{i_D(\fI_D)}
\egloz
determined by
\bgl
\label{inverse-piP}
  ((i_D)\dot{} \rte P) (\dot{v}_q^* i_{D / \fI_D}(d + \fI_D) \dot{v}_p) = v_q^* i_D(d) v_p + \spkl{i_D(\fI_D)}.
\egl
Comparing \eqref{piP} and \eqref{inverse-piP}, we see that $(i_D)\dot{} \rte P$ is the inverse of $(\pi \rte P)\dot{}$.
\eproof

\bsubcor
\label{pi3}
In the same situation as in the previous two lemmas, let $E: D \rte_{\alpha} P \to D$ be the faithful conditional expectation from Lemma~\ref{pi1}. Let $\fI$ be an ideal of $D \rte_{\alpha} P$ as in Lemma~\ref{pi2}. Then for all $x \in (D \rte_{\alpha} P)_{+}$, $E(x) \in \fI$ implies $x \in \fI$.
\esubcor
\bproof
Condition~\eqref{sat} is satisfied for $i_{D / \fI_D}$ and $\dot{v}$ from Lemma~\ref{pi2} since
\bgloz
  \dot{v}_p^* i_{D / \fI_D}(d + \fI_D) \dot{v}_p = (\pi \rte P)\dot{} \, (\underbrace{v_p^* i_D(d) v_p}_{
  \begin{array}{c}
  \in i_D(D) \\
  \text{by \eqref{sat}}
  \end{array}
  }
  )
  \in (\pi \rte P)\dot{} \, (i_D(D)) = i_{D / \fI_D}(D / \fI_D).
\egloz
Thus, by Lemma~\ref{pi1} applied to the C*-dynamical semisystem $(D / \fI_D,P,\dot{\alpha})$, there exists a faithful conditional expectation
\bgloz
  \dot{E}: (D / \fI_D) \rte_{\dot{\alpha}} P \to D / \fI_D
\egloz
which is determined by
\bgl
\label{dotE}
  \dot{E}(\dot{v_q}^* i_{D / \fI_D}(d + \fI_D) \dot{v}_p) = \delta_{q,p} (i_{D / \fI_D})^{-1} (\dot{v}_p^* i_{D / \fI_D}(d + \fI_D) \dot{v}_p).
\egl
Comparing \eqref{E-det} and \eqref{dotE} and using the homomorphism $\pi \rte P$ from Lemma~\ref{pi2}, we see that the diagram
\bgl
\label{E-dotE}
  \begin{CD}
  D \rte_{\alpha} P @> \pi \rte P >> (D / \fI_D) \rte_{\dot{\alpha}} P \\
  @V E VV @V \dot{E} VV \\
  D @> \pi >> D / \fI_D
  \end{CD}
\egl
commutes.

Now take $x \in (D \rte_{\alpha} P)_{+}$ with $E(x) \in \fI$. This means that
\bgloz
  0 + \fI_D = \pi(E(x)) \overset{\eqref{E-dotE}}{=} \dot{E}((\pi \rte P)(x)).
\egloz
As $\dot{E}$ is faithful, we conclude that $(\pi \rte P)(x) = 0$ in $(D / \fI_D) \rte_{\dot{\alpha}} P$. From Lemma~\ref{pi2}, it follows directly that the kernel of $\pi \rte P$ is $\spkl{i_D(\fI_D)}$. Thus we conclude that $x \in \spkl{i_D(\fI_D)} \subseteq \fI$.
\eproof

Let us now return to the situation of interest. Let $R$ be a Dedekind domain, and let $R \rtimes R\reg$ be the $ax+b$-semigroup over $R$. As explained in the second half of \S~2.4 in \cite{Li2}, the family of right ideals $\cJ$ of $R \rtimes R\reg$ is given by
\bgl
\label{cJax+b}
  \cJ = \menge{(r+I) \times I\reg}{r \in R, (0) \neq I \triangleleft R} \cup \gekl{\emptyset}.
\egl
By Proposition~3.13 in \cite{Li2} and because the constructible right ideals of $R \rtimes R\reg$ are independent by Lemma~2.29 in \cite{Li2}, the left regular representation
\bgloz
  \lambda: C^*(R \rtimes R\reg) \to C^*_r(R \rtimes R\reg)
\egloz
from \S~2.1 in \cite{Li2} is an isomorphism. Using $\lambda$, we will from now on always identify the full and the reduced semigroup C*-algebras of $R \rtimes R\reg$ (i.e. we may write $v_p$ for $V_p$ and $e_X$ for $E_X$ using the notation from \cite{Li2}). Moreover, for a subset $X$ of $R \rtimes R\reg$, let $e_X \in \cL(\ell^2(R \rtimes R\reg))$ be the orthogonal projection onto the subspace $\ell^2(X) \subseteq \ell^2(R \rtimes R\reg)$. This is consistent with the notation from \cite{Li2}. We sometimes write $\e{X}$ for $e_X$ if the expression for $X$ is rather long.

Now set
\bgloz
  D(R \rtimes R\reg) = C^*(\menge{e_{(r+I) \times I\reg}}{r \in R, (0) \neq I \triangleleft R}),
\egloz
and let $\tau$ be the action of $R \rtimes R\reg$ on $D(R \rtimes R\reg)$ given by $\tau_p = \Ad(v_p)$ for all $p \in R \rtimes R\reg$. By Lemma~2.14 in \cite{Li2}, we can canonically identify $C^*(R \rtimes R\reg)$ (hence also $C^*_r(R \rtimes R\reg)$) with $D(R \rtimes R\reg) \rte_{\tau} (R \rtimes R\reg)$. As $\tau_p$ is given by conjugation with an isometry, it is injective. Moreover, we have already seen in the proof of Theorem~\ref{Kax+b} that $R \rtimes R\reg$ is a left Ore semigroup whose enveloping group is amenable. As \eqref{sat} is also satisfied for $(D(R \rtimes R\reg),R \rtimes R\reg,\tau)$ by Corollary~2.9 in \cite{Li2}, we can apply Lemma~\ref{pi1} to $(D(R \rtimes R\reg),R \rtimes R\reg,\tau)$. Using the canonical identification of $C^*_r(R \rtimes R\reg)$ with $D(R \rtimes R\reg) \rte_{\tau} (R \rtimes R\reg)$, we obtain a faithful conditional expectation $E: C^*_r(R \rtimes R\reg) \to D(R \rtimes R\reg)$ which is determined by
\bgloz
  E(v_q^* e_{(r+I) \times I\reg} v_p) = \delta_{q,p} v_p^* e_{(r+I) \times I\reg} v_p.
\egloz
Corollary~\ref{pi3} then tells us the following:
\bsubcor
\label{pi3'}
Let $\fI$ be an ideal of $C^*_r(R \rtimes R\reg)$. For a positive element $y$ in $C^*_r(R \rtimes R\reg)$, $E(y) \in \fI$ implies $y \in \fI$.
\esubcor

Our first goal is to prove the following variation of Lemma~4.12 from \cite{C-D-L}:
\bsublemma
\label{pi4}
Let $R$ be a Dedekind domain with infinitely many pairwise distinct prime ideals, let $\fI$ be an ideal of $C^*_r(R \rtimes R\reg)$, and let $y$ be a positive element in the *-algebra generated by the isometries $v_p$, $p \in R \rtimes R\reg$, i.e.
\bgloz
  y \in (\stalg(\menge{v_p}{p \in R \rtimes R\reg}))_{+}.
\egloz

If $y$ does not lie in $\fI$, then there is a projection $\delta$ in $C^*_r(R \rtimes R\reg)$ of the form
\bgl
\label{delta}
  \delta = \e{(r+I) \times I\reg \setminus \bigcup_{k=1}^n (s_k + J_k) \times J_k\reg}
\egl
with $r, s_1, \dotsc, s_n \in R$ and non-zero ideals $I, J_1, \dotsc, J_n$ of $R$ such that
\begin{itemize}
\item[1.] $\delta$ does not lie in $\fI$
\item[2.] $\delta y \delta = \rukl{\norm{E(y) + \fI}_{C^*_r(R \rtimes R\reg) / \fI}} \delta$.
\end{itemize}

Note that in \eqref{delta}, the case $n=0$ is possible; it corresponds to $\delta = \e{(r+I) \times I\reg}$.
\esublemma
\bproof
As $y$ lies in $\stalg(\menge{v_p}{p \in R \rtimes R\reg})$, it is of the form
\bgl
\label{y}
  y = d + \sum_{i=1}^m v_{q_i}^* d_i v_{p_i}
\egl
with $p_1, \dotsc, p_m$, $q_1, \dotsc, q_m$ in $R \rtimes R\reg$ such that $q_i \neq p_i$ for all $1 \leq i \leq m$ and where $d$ and $d_i$ ($1 \leq i \leq m$) are finite linear combinations of the projections $e_X$, $X \in \cJ$. The condition $q_i \neq p_i$ for all $1 \leq i \leq m$ implies $E(y) = d$. Moreover, we can write $d$ as a finite sum $d = \sum_X \lambda_X e_X$. Now we can orthogonalize the projections $e_X$ which appear in this finite sum, and we obtain pairwise orthogonal projections $e_Y$ and a presentation $d = \sum_Y \mu_Y e_Y$. By Corollary~\ref{pi3'}, $y \notin \fI$ implies $d = E(y) \notin \fI$. Thus $0 < \norm{E(y) + \fI}_{C^*_r(R \rtimes R\reg) / \fI} = \sup(\gekl{\mu_Y})$ where the supremum is taken over all coefficients $\mu_Y$ corresponding to $e_Y \notin \fI$ appearing in the sum above which represents $d$. Since this sum is finite and the $\gekl{e_Y}$ are pairwise orthogonal, there exists a projection $e_Y$ such that the corresponding coefficient precisely coincides with $\norm{E(y) + \fI}_{C^*_r(R \rtimes R\reg) / \fI}$. This implies that this projection satisfies
\bgl
\label{ede}
  e_Y d e_Y = \rukl{\norm{E(y) + \fI}_{C^*_r(R \rtimes R\reg) / \fI}} e_Y.
\egl
Moreover, since the projections $\gekl{e_Y}$ were obtained by orthogonalizing the commuting projections $\gekl{e_X}$, the subset $Y$ of $R \rtimes R\reg$ must be of the form
\bgl
\label{Y}
  Y = (\ti{r} + \ti{I}) \times \ti{I}\reg \setminus \bigcup_{k=1}^n (\ti{s}_k + \ti{J}_k) \times \ti{J}_k\reg
\egl
with $\ti{r}, \ti{s}_1, \dotsc, \ti{s}_n \in R$ and non-zero ideals $\ti{I}, \ti{J}_1, \dotsc, \ti{J}_n$ of $R$ such that $\ti{J}_k \subseteq \ti{I}$ for all $1 \leq k \leq n$. The case $n=0$ is allowed; it corresponds to $Y = (\ti{r} + \ti{I}) \times \ti{I}\reg$. That $Y$ is of this form follows from \eqref{cJax+b}.

We now choose $(b,a) \in R \rtimes R\reg$ satisfying
\begin{itemize}
\item[1$_{b,a}$.] $v_{(b,a)} e_Y = e_Y v_{(b,a)}$
\item[2$_{b,a}$.] $v_{(b,a)} v_{(b,a)}^* v_{q_i}^* d_i v_{p_i} v_{(b,a)} v_{(b,a)}^* = 0$ for all $1 \leq i \leq m$.
\end{itemize}

Let $q_i = (b'_i,a'_i) \in R \rtimes R\reg$ and $p_i = (b_i,a_i) \in R \rtimes R\reg$. Then
\bglnoz
  && v_{(b,a)} v_{(b,a)}^* v_{q_i}^* d_i v_{p_i} v_{(b,a)} v_{(b,a)}^* \\
  &=& v_{q_i}^* v_{q_i} e_{(b+aR) \times (aR)\reg} v_{q_i}^* d_i v_{p_i} e_{(b+aR) \times (aR)\reg} v_{p_i}^* v_{p_i} \\
  &=& v_{q_i}^* d_i \e{(b'_i + a'_ib + a'_iaR) \times (a'_iaR)\reg} \e{(b_i + a_ib + a_iaR) \times (a_iaR)\reg} v_{p_i} \\
  &=& v_{q_i}^* d_i \e{\rukl{(b'_i + a'_ib + a'_iaR) \times (a'_iaR)\reg} \cap \rukl{(b_i + a_ib + a_iaR) \times (a_iaR)\reg}} v_{p_i}
\eglnoz
vanishes if
\bgl
\label{cond2}
  (b'_i + a'_ib + a'_iaR) \cap (b_i + a_ib + a_iaR) = \emptyset
\egl
for all $1 \leq i \leq m$. If for all $1 \leq i \leq m$, we have $b'_i + a'_ib - (b_i + a_ib) = (b'_i - b_i) + (a'_i - a_i)b \notin aR$, then certainly \eqref{cond2} holds.

We claim that we can choose $b \in R$ such that
\begin{itemize}
\item[1$_b$.] $b \in \ti{J}_1 \cap \dotsb \cap \ti{J}_n$ (or $b \in \ti{I}$ if $n=0$)
\item[2$_b$.] $(b'_i - b_i) + (a'_i - a_i)b \neq 0$ for all $1 \leq i \leq m$.
\end{itemize}
The reason is that we have by assumption $(b'_i,a'_i) \neq (b_i,a_i)$ for all $1 \leq i \leq m$. This implies that for all $1 \leq i \leq m$, either $a'_i = a_i \wedge b'_i \neq b_i$ or $a'_i \neq a_i$. If $a'_i = a_i \wedge b'_i \neq b_i$, then $(b'_i - b_i) + (a'_i - a_i)b = b'_i - b_i \neq 0$ for all $b \in R$, and if $a'_i \neq a_i$, then $(b'_i - b_i) + (a'_i - a_i)b \neq 0$ for all $b \in R$ with $b \neq -(a'_i - a_i)^{-1}(b'_i - b_i)$. This shows that there are only finitely many ring elements which do not satisfy 2$_b$. As against that, by our assumption that $R$ is a Dedekind domain with infinitely many pairwise distinct prime ideals, $\ti{J}_1 \cap \dotsb \cap \ti{J}_n$ (or $\ti{I}$ if $n=0$) is an infinite set. Thus we can find $b$ in $R$ satisfying 1$_b$. and 2$_b$. at the same time. Let us fix such a choice for $b \in R$.

As a next step, we claim that we can choose $a \in R\reg$ such that
\begin{itemize}
\item[1$_a$.] $a \in 1+ \ti{J}_1 \cap \dotsb \cap \ti{J}_n$ (or $a \in 1 + \ti{I}$ if $n=0$)
\item[2$_a$.] $(b'_i - b_i) + (a'_i - a_i)b \notin aR$ for all $1 \leq i \leq m$.
\end{itemize}
To see this, first note that if $\prod_{i=1}^m ((b'_i - b_i) + (a'_i - a_i)b)$ does not lie in $aR$, then 2$_a$ follows. By 2$_b$., the element $\prod_{i=1}^m ((b'_i - b_i) + (a'_i - a_i)b)$ is not zero. Thus, by our assumption that $R$ is a Dedekind domain with infinitely many pairwise distinct prime ideals, there exists a prime ideal $P$ of $R$ such that $\prod_{i=1}^m ((b'_i - b_i) + (a'_i - a_i)b)$ does not lie in $P$ and also $\ti{J}_1 \cap \dotsb \cap \ti{J}_n \nsubseteq P$ (or $\ti{I} \nsubseteq P$ if $n=0$). By the Chinese Remainder Theorem (see for example \cite{Neu}, Chapter~I, Theorem~(3.6)), there exists a non-zero element $a$ of the prime ideal $P$ such that $a \in 1 + \ti{J}_1 \cap \dotsb \cap \ti{J}_n$ (or $a \in 1 + \ti{I}$ if $n=0$). This $a$ obviously satisfies 1$_a$. and 2$_a$.

Finally, we claim that this choice for $(b,a) \in R \rtimes R\reg$ satisfies 1$_{b,a}$. and 2$_{b,a}$. First, by our observations, 2$_a$. implies \eqref{cond2}, hence 2$_{b,a}$. To prove 1$_{b,a}$., note that 1$_a$. implies that $aR$ is coprime to each of the ideals $\ti{I}, \ti{J}_1, \dotsc, \ti{J}_n$, so that
\bgln
\label{coprime}
  && (a \ti{r} + a \ti{I}) \times (a \ti{I})\reg \setminus \bigcup_{k=1}^n (a \ti{s}_k + a \ti{J}_k) \times (a \ti{J}_k)\reg \\
  &=& \rukl{(a \ti{r} + \ti{I}) \times \ti{I}\reg \setminus \bigcup_{k=1}^n (a \ti{s}_k + \ti{J}_k) \times \ti{J}_k\reg} \cap \rukl{aR \times (aR)\reg}. \nonumber
\egln
Thus
\bglnoz
  v_{(0,a)} e_Y v_{(0,a)}^* 
  &\overset{\eqref{Y}}{=}& \e{(a \ti{r} + a \ti{I}) \times (a \ti{I})\reg \setminus \bigcup_{k=1}^n (a \ti{s}_k + a \ti{J}_k) \times (a \ti{J}_k)\reg} \\
  &\overset{\eqref{coprime}}{=}& \e{(a \ti{r} + \ti{I}) \times \ti{I}\reg \setminus \bigcup_{k=1}^n (a \ti{s}_k + \ti{J}_k) \times \ti{J}_k\reg} e_{aR \times (aR)\reg} \\
  &=& \e{(a \ti{r} + \ti{I}) \times \ti{I}\reg \setminus \bigcup_{k=1}^n (a \ti{s}_k + \ti{J}_k) \times \ti{J}_k\reg} v_{(0,a)} v_{(0,a)}^* \\
  &\overset{\text{1}_a.}{=}& \e{(\ti{r} + \ti{I}) \times \ti{I}\reg \setminus \bigcup_{k=1}^n (\ti{s}_k + \ti{J}_k) \times \ti{J}_k\reg} v_{(0,a)} v_{(0,a)}^* \\
  &=& e_Y v_{(0,a)} v_{(0,a)}^*.
\eglnoz
Multiplication with $v_{(0,a)}$ from the right implies 
\bgloz
  v_{(0,a)} e_Y = e_Y v_{(0,a)}.
\egloz
Furthermore, 1$_b$ implies
\bgloz
  v_{(b,1)} e_Y = e_Y v_{(b,1)}.
\egloz
Thus we obtain
\bgloz
  v_{(b,a)} e_Y = v_{(b,1)} v_{(0,a)} e_Y = v_{(b,1)} e_Y v_{(0,a)} = e_Y v_{(b,1)} v_{(0,a)} = e_Y v_{(b,a)}.
\egloz
This proves 1$_{b,a}$. So we have seen that $(b,a)$ satisfies 1$_{b,a}$. and 2$_{b,a}$.

Now we set $\delta \defeq v_{(b,a)} e_Y v_{(b,a)}^*$. Then $\delta$ is a projection in $C^*_r(R \rtimes R\reg)$ of the desired form as in \eqref{delta} by construction. In addition, $\delta$ does not lie in $\fI$ because $e_Y$ does not lie in $\fI$. And finally, we have
\bgl
\label{deltavev}
  \delta = v_{(b,a)} e_Y v_{(b,a)}^* \overset{\text{1}_{b,a}.}{=} e_Y v_{(b,a)} v_{(b,a)}^* = v_{(b,a)} v_{(b,a)}^* e_Y
\egl
and hence
\bglnoz
  \delta y \delta &\overset{\eqref{y}}{=}& \delta d \delta + \sum_{i=1}^m \delta v_{q_i}^* d_i v_{p_i} \delta \\
  &\overset{\eqref{deltavev}}{=}& v_{(b,a)} v_{(b,a)}^* e_Y d e_Y v_{(b,a)} v_{(b,a)}^* + 
  \sum_{i=1}^m e_Y \underbrace{v_{(b,a)} v_{(b,a)}^* v_{q_i}^* d_i v_{p_i} v_{(b,a)} v_{(b,a)}^*}_{=0 \text{ by } \text{2}_{(b,a)}.} e_Y \\
  &\overset{\eqref{ede}}{=}& \rukl{\norm{E(y) + \fI}_{C^*_r(R \rtimes R\reg) / \fI}} e_Y v_{(b,a)} v_{(b,a)}^* \\
  &\overset{\eqref{deltavev}}{=}& \rukl{\norm{E(y) + \fI}_{C^*_r(R \rtimes R\reg) / \fI}} \delta.
\eglnoz
Therefore this projection $\delta$ has the desired properties and satisfies conditions 1. and 2. of the lemma.
\eproof

To proceed, we need
\bsublemma
\label{pi5}
In the situation of the previous lemmma, the projection $\delta$ gives rise to a properly infinite projection $\delta + \fI$ in the quotient $C^*_r(R \rtimes R\reg) / \fI$.
\esublemma
\bproof
The projection $\delta$ is of the form $\delta = \e{(r+I) \times I\reg \setminus \bigcup_{k=1}^n (s_k + J_k) \times J_k\reg}$ by \eqref{delta}. Again, $n=0$ is allowed. Now choose $c \in R\reg$ and $r_1,r_2 \in R$ such that
\begin{itemize}
\item[(*$_c$)] $c$ is not invertible and $c$ lies in $1 + \bigcap_{k=1}^n J_k$ (or $1 + I$ if $n=0$)
\item[(*$_r$)] $r_1$, $r_2$ lie in $\bigcap_{k=1}^n J_k$ (or $I$ if $n=0$) but $r_1 + cR \neq r_2 + cR$.
\end{itemize}
This is possible because by our assumption that $R$ is a Dedekind domain with infinitely many pairwise distinct prime ideals, we can first find an element $c \in R\reg$ satisfying (*$_c$) using strong approximation (compare \cite{Bour}, Chapitre~VII, \S~2.4, Proposition~2). Then, as (*$_c$) implies that $cR$ and $\bigcap_{k=1}^n J_k$ (or $I$ for $n=0$) are coprime, the Chinese Remainder Theorem (see for instance \cite{Neu}, Chapter~I, Theorem~(3.6)) tells us that we can find elements $r_1$ and $r_2$ in $R$ satisfying (*$_r$). Then, by analogous computations as in the proof of the previous lemma, (*$_c$) and (*$_r$) imply
\bgl
\label{vd=dv}
  v_{(r_i,c)} \delta = \delta v_{(r_i,c)} \text{ for } i=1,2.
\egl

Set $\delta_i = v_{(r_i,c)} \delta v_{(r_i,c)}^*$ for $i=1,2$. Then we certainly have $\delta_i \sim \delta$ for $i=1,2$, where $\sim$ stands for \an{Murray-von Neumann equivalent}. Moreover, for $i=1,2$, we obtain
\bgloz
  \delta_i = v_{(r_i,c)} \delta v_{(r_i,c)}^* \overset{\eqref{vd=dv}}{=} \delta v_{(r_i,c)} v_{(r_i,c)}^* \leq \delta.
\egloz
And finally,
\bglnoz
  \delta_1 \delta_2 &=& v_{(r_1,c)} \delta v_{(r_1,d)}^* v_{(r_2,c)} \delta v_{(r_2,c)}^* \\
  &\overset{\eqref{vd=dv}}{=}& \delta v_{(r_1,c)} v_{(r_1,d)}^* v_{(r_2,c)} v_{(r_2,c)}^* \delta \\
  &=& \delta e_{(r_1+cR) \times (cR)\reg} e_{(r_2+cR) \times (cR)\reg} \delta \\
  &=& \delta \e{\rukl{(r_1+cR) \times (cR)\reg} \cap \rukl{(r_2+cR) \times (cR)\reg}} \delta \\
  &=& 0
\eglnoz
as $r_1 + cR \neq r_2 + cR$ by (*$_r$). As $\delta + \fI$ is a non-zero projection in the quotient $C^*_r(R \rtimes R\reg) / \fI$ by condition 1. in Lemma~\ref{pi4}, this proves our claim.
\eproof

With these preparations, we are ready for the
\bproof[Proof of Theorem~\ref{pi}]
Let $R$ be a Dedekind domain with infinitely many pairwise distinct prime ideals. By \cite{Pas-Ror}, Proposition~2.11, we have to prove that every non-zero hereditary sub-C*-algebra in any quotient of $C^*_r(R \rtimes R\reg)$ contains an infinite projection. Let $\fI$ be an ideal of $C^*_r(R \rtimes R\reg)$. It suffices to show that every hereditary sub-C*-algebra of $C^*_r(R \rtimes R\reg) / \fI$ of the form
\bgloz
  \overline{(z + \fI) \rukl{C^*_r(R \rtimes R\reg) / \fI} (z + \fI)}
\egloz
for some $z \in C^*_r(R \rtimes R\reg)_{+}$, $z \notin \fI$ contains an infinite projection because every non-zero hereditary sub-C*-algebra of $C^*_r(R \rtimes R\reg) / \fI$ contains a subalgebra of such a form.

First of all, $z \notin \fI$ implies $E(z) \notin \fI$ by Corollary~\ref{pi3'}. As $\stalg(\menge{v_p}{p \in R \rtimes R\reg})$ is dense in $C^*_r(R \rtimes R\reg)$, there exists a positive element $y$ in $\stalg(\menge{v_p}{p \in R \rtimes R\reg})$ such that
\bgloz
  \norm{z-y} < \tfrac{1}{3} \norm{E(z) + \fI}_{C^*_r(R \rtimes R\reg) / \fI}.
\egloz
It follows that
\bglnoz
  && \norm{E(z) - E(y) + \fI}_{C^*_r(R \rtimes R\reg) / \fI} \\
  &\leq& \norm{E(z) - E(y)} \leq \norm{z-y} \\
  &<& \tfrac{1}{3} \norm{E(z) + \fI}_{C^*_r(R \rtimes R\reg) / \fI}
\eglnoz
so that
\bgl
\label{Ey-ineq}
  \norm{E(y) + \fI}_{C^*_r(R \rtimes R\reg) / \fI} > \tfrac{2}{3} \norm{E(z) + \fI}_{C^*_r(R \rtimes R\reg) / \fI} > 0.
\egl
This implies that $E(y)$ does not lie in $\fI$ and hence, by Corollary~\ref{pi3'}, $y$ does not lie in $\fI$. By Lemma~\ref{pi4} and Lemma~\ref{pi5}, there exists a projection $\delta$ in $C^*_r(R \rtimes R\reg)$ such that $\delta + \fI$ is a properly infinite projection in $C^*_r(R \rtimes R\reg) / \fI$ and
\bgl
\label{dyd}
  \delta y \delta = \rukl{\norm{E(y) + \fI}_{C^*_r(R \rtimes R\reg) / \fI}} \delta.
\egl
We get
\bgloz
  \norm{\delta z \delta - \delta y \delta + \fI}_{C^*_r(R \rtimes R\reg) / \fI} \leq \norm{\delta} \norm{z-y} \norm{\delta} 
  < \tfrac{1}{3} \norm{E(z) + \fI}_{C^*_r(R \rtimes R\reg) / \fI}
\egloz
and
\bglnoz
  && \rukl{\delta y \delta - \tfrac{1}{3} \norm{E(z) + \fI}_{C^*_r(R \rtimes R\reg) / \fI}}_{+} \\
  &\overset{\eqref{dyd}}{=}& \rukl{\rukl{\norm{E(y) + \fI}_{C^*_r(R \rtimes R\reg) / \fI}} \delta - \tfrac{1}{3} \norm{E(z) + \fI}_{C^*_r(R \rtimes R\reg) / \fI}}_{+} \\
  &=& \underbrace{\rukl{\norm{E(y) + \fI}_{C^*_r(R \rtimes R\reg) / \fI} - \tfrac{1}{3} \norm{E(z) + \fI}_{C^*_r(R \rtimes R\reg) / \fI}}}_{\eqdef C} \delta
\eglnoz
with $C > \tfrac{1}{3} \norm{E(z) + \fI}_{C^*_r(R \rtimes R\reg) / \fI} > 0$ by \eqref{Ey-ineq}. In this situation, by Lemma~2.2 in \cite{Kir-Ror2} (applied to $A = C^*_r(R \rtimes R\reg) / \fI$, $a = \delta y \delta + \fI$, $b = \delta z \delta + \fI$ and $\ve = \tfrac{1}{3} \norm{E(z) + \fI}_{C^*_r(R \rtimes R\reg) / \fI}$), there exists $x' \in C^*_r(R \rtimes R\reg) / \fI$ with $C \delta + \fI = x' (\delta z \delta + \fI) x'^*$. Now set $x \defeq C^{-\tfrac{1}{2}} x' (\delta + \fI)$. Then $\delta + \fI = x (z + \fI) x^*$ is a properly infinite projection (see Lemma~\ref{pi5}). We conclude that
\bgloz
  (z + \fI)^{\tfrac{1}{2}} x^* x (z + \fI)^{\tfrac{1}{2}}
\egloz
is a projection in
\bgloz
  \overline{(z + \fI) \rukl{C^*_r(R \rtimes R\reg) / \fI} (z + \fI)}
\egloz
which is Murray-von Neumann equivalent to $\delta + \fI$, hence properly infinite itself.
\eproof

Combining Theorem~\ref{Kax+b} and Theorem~\ref{pi} with the K-theoretic results for ring C*-algebras from \cite{Cu-Li} and \cite{Li-Lu}, we obtain
\bsubcor
For every ring of integers $R$ in a number field, the semigroup C*-algebra $C^*_r(R \rtimes R\reg)$ is purely infinite, has the ideal property but does not have real rank zero.
\esubcor
\bproof
Let $R$ be the ring of integers in a number field. Comparing universal properties, it is clear that the ring C*-algebra $\fA[R]$ of $R$ is a quotient of the semigroup C*-algebra $C^*(R \rtimes R\reg)$ of the $ax+b$-semigroup over $R$. Thus $\fA[R]$ is also a quotient of $C^*_r(R \rtimes R\reg)$. We have proven in \cite{Cu-Li} and \cite{Li-Lu} that $K_0(\fA[R])$ cannot be finitely generated, whereas it follows from Theorem~\ref{Kax+b} that $K_0(C^*_r(R \rtimes R\reg))$ is finitely generated. Hence the quotient map from $C^*_r(R \rtimes R\reg)$ to $\fA[R]$ cannot be surjective on $K_0$. In the language of \cite{Pas-Ror}, this means that $C^*_r(R \rtimes R\reg)$ is not $K_0$-liftable. As we have seen in Theorem~\ref{pi} that $C^*_r(R \rtimes R\reg)$ is purely infinite, Theorem~4.2 in \cite{Pas-Ror} implies that $C^*_r(R \rtimes R\reg)$ cannot have real rank zero. This shows the last part of our assertion. The first part follows from Theorem~\ref{pi}.
\eproof

{\bf Remark}
\\
For a cancellative semigroup, we do not only have the left regular representation, but also the right regular one. For groups, the C*-algebras generated by these representations are isomorphic due to invertibility of the group elements. But for a genuine (and let us say non-abelian) semigroup, the left and right regular representations generate in general different C*-algebras. For our present piece of work, the C*-algebra generated by the left regular representation of the $ax+b$-semigroup over the ring of integers of a number field was the motivating example. A natural question would be:

What about the C*-algebra generated by the right regular representation of such an $ax+b$-semigroup?

It turns out that although the C*-algebras for the left and right regular repesentations of such semigroups are quite different (the one for the right regular representation is not purely infinite), their K-theoretic invariants do coincide. In a forthcoming paper, the authors plan to discuss the C*-algebras of the right regular representations of such $ax+b$-semigroups in a general context.

\end{document}